\documentclass{IEEEtran}
\usepackage{amsmath}
\usepackage{amssymb}
\usepackage{array}
\usepackage{multirow}
\usepackage{latexsym}
\usepackage{algorithm, algpseudocode}
\usepackage{graphicx}
\usepackage[caption=false,subrefformat=parens,labelformat=parens,font=footnotesize]{subfig}
\usepackage{color}
\usepackage{url}
\usepackage{fancyvrb}
\usepackage{rotating}
\usepackage{mathtools}
\usepackage{xspace}
\usepackage{longtable}

\algrenewcommand{\algorithmiccomment}[1]{\hfill$\blacktriangleright$ #1}

\def\pdffigdir{figures-pdf}

\ifdefined\tikzdir
  \usepackage{pgfplots,tikzscale}
  \usetikzlibrary{external}
  \usetikzlibrary{arrows}
  \usetikzlibrary{calc}
  \usetikzlibrary{intersections}
  \pgfplotsset{compat=newest}
  \pgfplotsset{plot coordinates/math parser=false}
  \newcommand{\includepgfplots}[2][width=1\textwidth]{%
    \tikzsetnextfilename{#2}%
    \includegraphics[#1]{\tikzdir/#2.tikz}}
\else
  \newcommand{\includepgfplots}[2][]{%
    \includegraphics{\pdffigdir/#2}}
\fi

\newcolumntype{L}[1]{>{\raggedright\let\newline\\\arraybackslash\hspace{0pt}}m{#1}}
\newcolumntype{C}[1]{>{\centering\let\newline\\\arraybackslash\hspace{0pt}}m{#1}}
\newcolumntype{R}[1]{>{\raggedleft\let\newline\\\arraybackslash\hspace{0pt}}m{#1}}

\newcommand{\beq}{\begin{equation}}
\newcommand{\eeq}{\end{equation}}
\newcommand{\bseq}{\begin{subequations}}
\newcommand{\eseq}{\end{subequations}}

\newcommand{\be}{\begin{enumerate}}
\newcommand{\ee}{\end{enumerate}}
\newcommand{\bi}{\begin{itemize}}
\newcommand{\ei}{\end{itemize}}

\newcommand{\bR}{{\mathbb R}}

\newcommand{\cA}{{\cal A}}

\newcommand{\cD}{{\cal D}}

\newcommand{\cF}{{\cal F}}
\newcommand{\cG}{{\cal G}}
\newcommand{\cH}{{\cal H}}

\newcommand{\cL}{{\cal L}}

\newcommand{\cN}{{\cal N}}

\newcommand{\cS}{{\cal S}}
\newcommand{\cT}{{\cal T}}

\newcommand{\cW}{{\cal W}}

\DeclareMathOperator*{\s.t.}{subject\ to}

\DeclareMathOperator*{\argmax}{arg\,max}

\newcommand{\lb}[1]{\underline{#1}}
\newcommand{\ub}[1]{\overline{#1}}

\newcommand{\omin}[1]{\min_{#1} \quad}
\newcommand{\omax}[1]{\max_{#1} \quad}

\newcommand{\oargmax}[1]{\argmax_{#1} \quad}
\newcommand{\ost}{\s.t. \quad}

\newcommand{\bm}[1]{\begin{bmatrix}#1\end{bmatrix}}

\def\<{\left<}
\def\>{\right>}
\def\({\left(}
\def\){\right)}

\def\mc{\multicolumn}

\def\sjwresolved#1{}

\def\noprint#1{}

\newcommand{\tcell}[2][c]{\begin{tabular}[#1]{@{}c@{}}#2\end{tabular}}

\def\MATPOWER{\textsc{MatPower}\xspace}
\def\MATLAB{\textsc{Matlab}\xspace}
\def\IPOPT{\textsc{Ipopt}\xspace}

\everymath{\displaystyle}
\begin{document}

\title{Vulnerability Analysis of Power Systems\footnote{\today}}


\author{Taedong~Kim\IEEEauthorrefmark{2},
  Stephen~J.~Wright\IEEEauthorrefmark{2},
  Daniel~Bienstock\IEEEauthorrefmark{3}, and
  Sean~Harnett\IEEEauthorrefmark{4}%
\thanks{This work is supported by
  a DOE grant subcontracted through Argonne National Laboratory Award
  3F-30222, and National Science Foundation Grant DMS-1216318.}}

\maketitle

\renewcommand{\thefootnote}{\fnsymbol{footnote}}

\footnotetext[2]{Computer Sciences Department, 1210
  W. Dayton Street, University of Wisconsin, Madison, WI 53706,
  USA (email: \mbox{tdkim@cs.wisc.edu}; \mbox{swright@cs.wisc.edu})}
\footnotetext[3]{Department of Industrial Engineering and Operations Research
  and Department of Applied Physics and Applied Mathematics, Columbia University,
  500 West 120th St. New York, NY 10027, USA (email: \mbox{dano@columbia.edu})}
\footnotetext[4]{Department of Applied Physics and Applied Mathematics,
  Columbia University, 500 West 120th St. New York, NY 10027, USA
  (email: \mbox{sharnett@gmail.com})}

\renewcommand{\thefootnote}{\arabic{footnote}}

\begin{abstract}
Potential vulnerabilities in a power grid can be exposed by
identifying those transmission lines on which attacks (in the form of
interference with their transmission capabilities) causes maximum
disruption to the grid. In this study, we model the grid by
(nonlinear) AC power flow equations, and assume that attacks take the
form of increased impedance along transmission lines. We quantify
disruption in several different ways, including (a) overall deviation
of the voltages at the buses from $1.0$ per unit (p.u.), and (b) the
minimal amount of load that must be shed in order to restore the grid
to stable operation. We describe optimization formulations of the
problem of finding the most disruptive attack, which are either
nonlinear programing problems or nonlinear bilevel optimization
problems, and describe customized algorithms for solving these
problems.  Experimental results on the IEEE 118-Bus system and a
Polish 2383-Bus system are presented.
\end{abstract}


\section{Introduction}\label{sec:intro}

Identifying the vulnerable components in a power grid is vital to the
design and operation of a secure, stable system. One aspect of
vulnerability analysis is to identify those transmission lines for
which minor perturbations in their conductive properties leads to
major disruptions to the grid, such as voltage drops, or the need for
load shedding at demand nodes to restore feasible operation.

Vulnerability assessment for power systems has been widely studied in
recent times. Most works focus on minimizing the costs of load
shedding and additional generation in the DC model (which is
relatively easy to solve) or in lossless AC models (still relatively
easy to solve and analyze). In \cite{SalWB04,MotAG05,ArrG05},
identification of critical components of a power system is formulated
in a mixed-integer bilevel programming framework, and attacks on
different types of system components (transmission lines, generators,
and transformers) are considered.
The lower-level problem in the bilevel formulation is replaced by its
dual in \cite{MotAG05} and is approximated using KKT conditions in
\cite{ArrG05}. As an extension of \cite{SalWB04}, an approach based on
Bender's decomposition is proposed to solve larger instances of the
transmission line attack problem in \cite{SalWB09}.

Vulnerability assessment using the lossless AC model is studied in
\cite{DonLL08,PinMD10}. In these papers, transmission-line attacks are
formulated as bilevel optimization problems, in which either unmet
demands are maximized or attack costs (number of lines to attack) are
minimized to meet a specified level of grid disruption. (These models
are also discussed in \cite{Arr10}, which describes the equivalence of
the two models.) Then the lower-level problem is replaced by its KKT
conditions, yielding the single-level optimization problem that is
actually solved. In \cite{DonLL08}, this mixed-integer problem is
relaxed to a continuous problem (the binary variables are relaxed to
real variables confined to the interval $[0,1]$), while \cite{PinMD10}
develops a graph-partitioning approach by identifying load-rich and
generation-rich regions.

The paper \cite{DelAA10} describes a model that uses both load shedding
and line switching as defensive operations to reduce the disruption of
the system; the model is solved via Bender's decomposition with a
restart framework. Use of a genetic algorithm to solve the ``$N-k$''
problem (identifying the set of $k$ lines in a grid of $N$ lines whose
removal causes maximum disruption) is discussed in \cite{ArrF13}. A
minimum-cardinality approach (solved using a cutting-plane method) and
a continuous nonlinear attack model employing the DC power flow to
represent power grids, where a fictitious adversary modifies
reactances, are applied to the ``$N-k$'' problem in \cite{BieV10}.

In this paper, we propose two optimization models for vulnerability
analysis. Both models are founded on the AC power flow equations, and
both consider attacks in which the impedances of transmission lines
are increased.  In both formulations, the attacks respect a certain
``budget;'' their total amount of impedance adjustment cannot exceed a
certain specified level. The goal of the attacks is to maximize
disruption, as measured by two different metrics.

The first metric quantifies voltage disturbance at the buses, leading
to a nonlinear programming formulation. The voltage disturbance
usually appears as voltage drop, which often leads to an undesirable
situation where voltages become low enough that the system cannot
maintain stability. This situation, which is called voltage collapse
or voltage instability, can happen either quickly or relatively
slowly, and is characterized by a parallel process where reactive
power demand correspondingly increases. This eventuality causes the
active-power behavior of the system to approach the "nose" of the
$P-V$ curve. A more complete description is provided in \cite[pages 31
  and 35]{usc04}. With our first metric, we estimate this possible
voltage instability of a power grid assuming there is no response from
a system operator to the attack.

The second metric we consider here is a weighted sum of the amount of
load shedding (at demand nodes) and generation reduction (at
generation nodes) that is required to restore feasible operation of
the grid following the attack. This power adjustment is considered as
a defensive action of a system operator to keep voltages within a
stable range to avoid voltage collapse. This case is modeled as a
bilevel optimization problem, in which the lower level finds the
minimum load adjustment required to respond to the attack, and the
upper-level problem is to find the most disruptive attack.

In some existing literature, including some of the papers cited above,
a bilevel optimization model is reformulated as a single-level
optimization problem by replacing the lower level problem by its
optimality conditions. This formulation strategy is unappealing, as
the optimality conditions characterize only a stationary point, rather
than a minimizer, so they may allow consideration of saddle points or
local maximizers. In addition, if the bilevel formulation is designed
to solve the attacker-defender framework that we consider in this
paper, the reformulated single-level optimization model constructed by
replacing the lower-level problem by primal-dual optimality conditions
has the serious flaw that the model may exclude the most effective
attack. Specifically, an attack (upper level decision) that leads to
an infeasible lower level problem obviously maximizes the disruption
and thus is ``optimal'' for the attack problem (since it is not
possible to make a operational decision at the lower level to defend
against the attack). However, such an attack is excluded from
consideration by the single-level reformulation since no (lower-level)
primal-dual point satisfies the optimality condition constraints under
the attack. Thus, the single-level formulation will ignore the most
critical attack. Another drawback of single-level reformulations is
that the optimality-condition constraints may violate constraint
qualifications, causing possible complications in convergence
behavior.

The main contributions of this paper can be summarized as follows:
\be
 \item In contrast to previous attack models, the grid is modeled with
   full AC power flow equations, which are the most accurate
   mathematical models of power flow.
\item In our bilevel optimization formulation, we actually {\em solve}
  the lower-level problem rather than replacing it by its optimality
  conditions, as is done in earlier works, to avoid the formulation
  defects discussed above.
\item We develop effective heuristics that make our formulations
  tractable even for power grids with thousands of buses.
\ee

The remaining sections are organized as follows. We develop the
optimization models in Section~\ref{sec:problem} and describe the
challenges to be addressed in solving them. Section~\ref{sec:alg}
describes heuristics and optimization techniques that address these
challenges and that yield solutions of the problems. Experimental
results on 118-Bus and 2383-Bus cases are presented in
Section~\ref{sec:result}, and we discuss conclusions in
Section~\ref{sec:conclusion}.

\section{Problem Description}\label{sec:problem}

In this section we discuss power systems background and notation, and
describe our two formulations of the vulnerability analysis problem.
We describe notation and background on power flow equations in
Subsection~\ref{sec:notation}. Our first vulnerability model, based on
a voltage disturbance objective, is discussed in
Subsection~\ref{sec:problem.volt}. The second model, based on a
power-adjustment criterion, is discussed in
Subsection~\ref{sec:problem.load}.


\subsection{Notations and Background}\label{sec:notation}

We summarize here the power systems notation used in later sections,
most of which is standard.
\bi
  \setlength{\itemindent}{-1em}
  \item Set of buses: $\cN$
  \item Set of generators: $\cG \subseteq \cN$
  \item Set of demand buses: $\cD \subseteq \cN$
  \item Index of the slack bus: $s\in\cN$
  \item Set of transmission lines: $\cL \subseteq \cN\times\cN$
  \item Unit imaginary number: $j= \sqrt{-1}$
  \item Complex power at bus $i\in\cN$: $P_i+jQ_i$ (active power: $P_i$; real power: $Q_i$)
  \item Complex voltage at bus $i\in\cN$: $V_ie^{j\theta_i}$ (voltage magnitude: $V_i$; phase angle: $\theta_i$)
  \item Difference of angles $\theta_{i}$ and $\theta_{i'}$, for $(i,i')\in \cL$: $\theta_{ii'}:=\theta_{i}-\theta_{i'}$
  \item $(i,i')$ entry of the admittance matrix for the (unperturbed)
    grid: $G_{ii'}+jB_{ii'}$ (conductance: $G_{ii'}$; susceptance:
    $B_{ii'}$).
\ei
We assume that the set of generators $\cG$ and the set of demand buses
$\cD$ form a partition of $\cN$.

An attack on the grid is specified by means of a line perturbation
vector: $\gamma\in\bR_+^{|\cL|}$, with $\gamma_{ii'}$ denoting the
relative increase in impedance on line $(i,i')\in\cL$. Specifically,
an attack designated by the vector $\gamma$ causes conductances and
susceptances to be modified as follows:
\begin{align*}
G_{ii'}(\gamma) &=
  \begin{cases}
    \frac{G_{ii'}}{\gamma_{ii'}+1} & \text{if }i\neq i', \\
    -\sum_{i\neq m}G_{im}(\gamma)  & \text{if }i=i',
  \end{cases}\\
B_{ii'}(\gamma) &=
  \begin{cases}
    \frac{B_{ii'}}{\gamma_{ii'}+1}                                      & \text{if }i\neq i', \\
    -\sum_{i\neq m}\left(B_{im}(\gamma) - \frac{1}{2}B_{im}^{sh}\right) & \text{if }i=i',
  \end{cases}
\end{align*}
where $B_{im}^{sh}$ is the shunt (line charging) admittance of line
$(i,m)\in\cL$. More details on the bus admittance matrix can be found
from \cite[Chapter 9]{BerV99}. Note that when $\gamma_{ii'}=0$ for all
$(i,i')\in\cL$, the conductances and susceptances all attain their
original (unperturbed) values.


The AC power flow equations with perturbations $\gamma$ can be written
as follows:
\begin{equation}
  \bm{ F^P(V,\theta;\gamma)\\ F^Q(V,\theta;\gamma) } \, = \, 0,
\end{equation}
where the entries of $F^P$ and $F^Q$ (for all $i \in \cN$) are defined
as follows, for all $i \in \cN$:
\bseq\label{eq:power.equation}
\begin{align}
\begin{split}
 F^P_i(&V,\theta;\gamma) \, := \\
    & V_i\sum_{\mathclap{i':(i,i')\in\cL}}V_{i'} (G_{ii'}(\gamma)\cos{\theta_{ii'}}+B_{ii'}(\gamma)\sin{\theta_{ii'}})-P_i,
\end{split} \\
\begin{split}
\label{eq:power.equation.Q}
 F^Q_i(&V,\theta;\gamma) \, := \\
    & V_i\sum_{\mathclap{i':(i,i')\in\cL}}V_{i'} (G_{ii'}(\gamma)\sin{\theta_{ii'}}-B_{ii'}(\gamma)\cos{\theta_{ii'}})-Q_i.
\end{split}
\end{align}
\eseq
We assume throughout the paper that $P_i>0$ for generator buses
$i\in\cG$ and $P_i<0$ and $Q_i<0$ for demand buses $i\in\cD$. The
power flow problem is to find the values of the vectors $V$, $\theta$,
$P$, and $Q$ that satisfy equations \eqref{eq:power.equation}, given
the load demands $P_\cD$ and $Q_\cD$ at load buses and the voltage
magnitudes $V_\cG$ and active power injection $P_\cG$ at the generator
buses.  Conventionally, the reactive powers $Q_{\cG}$ are eliminated
from the problem (since they can be obtained explicitly from
\eqref{eq:power.equation.Q} for $i \in \cG$, and appear in no other
equations), yielding the following  reduced formulation:
\beq\label{eq:power.flow}
F(V,\theta;\gamma) \, = \,
\bm{ F_\cG^P(V,\theta;\gamma)\\
     F_\cD^P(V,\theta;\gamma)\\
     F_\cD^Q(V,\theta;\gamma)}
  \, = \, 0.
\eeq
Here, $V_s$, $\theta_s$, $V_{\cG}$, $P_{\cG}$, $P_{\cD}$, and $Q_{\cD}$
are parameters associated with the network; $\gamma$ is the impedance
modification vector described above; and $V_{\cD}$ and
$\theta_{\cG\cup\cD}$ are the variables in the model. These equations
usually can be solved using Newton's method, when the system has a
solution.
For additional details of formulation of power flow problems, see
\cite[Chapter 10]{BerV99}.

\subsection{Voltage Disturbance Model}\label{sec:problem.volt}

The AC power flow problem \eqref{eq:power.flow} often has multiple
solutions \cite{IbaIT80}, but only those solutions with $V_i\approx
1.0$ per unit (p.u.) for all $i \in \cD$ are stable and operational in
practice. In the vulnerability model described in this subsection, we
use the sum-of-squares deviation $\cF_V$ of the voltages from $1.0$
p.u.  as a measure of the disruption caused by an attack:
\begin{equation} \label{eq:def.FV}
\cF_V(\gamma) :=
\begin{cases}
  \frac{1}{2}\sum_{i\in\cD}(V_i-1)^2 & \parbox{3.5cm}{where $V$ is obtained by solving $F(V,\theta;\gamma)=0$,}\\
  +\infty            & \parbox{3.5cm}{when $F(V,\theta;\gamma)=0$ has no solution.}
\end{cases}
\end{equation}
Here, $\cF_V$ is a function of $\gamma$, the vector of relative
impedance increases. Note that only the voltage magnitudes of demand
buses $\cD$ are considered in $\cF_V(\gamma)$, since the voltage
magnitudes for generators and slack bus are given and fixed. We set
$\cF_V(\gamma)=+\infty$ when the attack results in an infeasible grid,
since such attacks are the best possible.

To limit the power of the purported attacker, we impose a constraint
on the vector $\gamma$, and define the {\em voltage disturbance
  vulnerability problem} as follows:
\bseq \label{eq:model.volt}
\begin{align}
\cH_V(\kappa, \ub{\gamma}) :=
\omax{\gamma} & \cF_V(\gamma) \\
\ost & e^T\gamma \le \kappa\ub{\gamma} \label{eq:model.volt.budget} \\
     & 0 \le \gamma \le \ub{\gamma}e, \label{eq:model.volt.limit}
\end{align}
\eseq
where $e=(1,1,\dotsc,1)^T$, the scalar $\ub{\gamma}$ is an upper bound
on relative impedance perturbation for each line, and $\kappa$ is the
maximum number of lines that can be attacked at the maximum
level. (Note that the actual number of lines attacked may be greater
than $\kappa$ if non-maximal attacks are made on some lines.)

Note that although the following model is a plausible alternative to
\eqref{eq:model.volt}, it is actually not valid:
\bseq\label{eq:model.volt.nlp}
\begin{alignat}{1}
\omax{V_\cD, \theta_{\cD\cup\cG}, \gamma} & \frac{1}{2}\sum_{i\in\cD}(V_i-1)^2 \\
\ost & F(V,\theta;\gamma) = 0 \\
    & e^T\gamma \le \kappa\ub{\gamma} \\
    & 0 \le \gamma \le \ub{\gamma}e.
\end{alignat}
\eseq
The reason is that when there is an attack $\gamma$ satisfying
\eqref{eq:model.volt.budget} and \eqref{eq:model.volt.limit} that
results in an infeasible grid, the formulation \eqref{eq:model.volt}
will find it (with an objective function of $+\infty$) while the
formulation \eqref{eq:model.volt.nlp} will not. In other words, the
formulation \eqref{eq:model.volt.nlp} does not fully capture the
adversarial nature of the attack. However, as a practical matter,
these two formulations find the same solution in cases where every
$\gamma$ satisfying \eqref{eq:model.volt.budget} and
\eqref{eq:model.volt.limit} allows for a feasible solution of the AC
power flow equations.

\subsection{Power-Adjustment Model}
\label{sec:problem.load}

Our second way to measure severity of an attack is to consider the
minimum adjustments to power that must be made to restore the grid to
feasible operation. Power adjustments take the form of shedding load
at demand nodes and adjusting generation at generator nodes. (We use
weights in the objective to discourage adjustment on nodes where it is
undesirable, such as at generators whose output cannot be adjusted or
at critical demand nodes whose load cannot be changed.)  Calculation
of this weighted sum of power adjustments involves solving a nonlinear
programming problem that we call the {\em feasibility restoration}
problem. This problem forms the lower-level problem in the bilevel
optimization problem, as we outline at the end of this subsection.

\subsubsection{Feasibility Restoration}\label{sec:problem.fr}

When the attack represented by $\gamma$ is too severe, the AC power
flow equations \eqref{eq:power.flow} may not have a solution for which
the voltages lie within an acceptable range. The {\em feasibility
  restoration} problem finds minimal adjustments to the power demands
(at demand nodes $\cD$) and power generation (at generator nodes
$\cG$) for which feasibility is restored to the AC power flow
equations. The formulation is as follows:
\bseq \label{eq:feas.rest.power}
\begin{alignat}{2}
\cF_L(\gamma) := \omin{\mathclap{\substack{V_\cD,\theta_{\cD\cup\cG},\\\sigma_\cG^+, \sigma_\cG^-,\rho_\cD}}} &
   \mathrlap{\sum_{i\in \cG} {\omega_i |P_i|(\sigma_i^++\sigma_i^-)}
              + \sum_{i\in \cD} {\omega_i |P_i|\rho_i}}\label{eq:feas.rest.power.obj}\\
\ost & F^P_\cG(V,\theta;\gamma)-{|P_\cG|\odot(\sigma_\cG^+-\sigma_\cG^-)} = 0  \label{eq:feas.rest.power.pv.p}\\
   & F^P_\cD(V,\theta;\gamma)-{|P_\cD|\odot\rho_\cD} = 0  \label{eq:feas.rest.power.pq.p}\\
   & F^Q_\cD(V,\theta;\gamma)-{|Q_\cD|\odot\rho_\cD} = 0  \label{eq:feas.rest.power.pq.q}\\
   & \lb{V} \le V_\cD \le \ub{V}                     \label{eq:feas.rest.power.v}\\
   & 0 \le {\sigma_\cG^+ \le \ub{\sigma}_\cG^+}      \label{eq:feas.rest.power.sigma.p}\\
   & 0 \le {\sigma_\cG^- \le \ub{\sigma}_\cG^-}      \label{eq:feas.rest.power.sigma.m}\\
   & 0 \le \rho_\cD \le \ub{\rho}_\cD,               \label{eq:feas.rest.power.rho}
\end{alignat}
\eseq
where $a\odot b$ is element-wise multiplication of vectors $a$ and
$b$.
Here, the variables $\sigma^+, \sigma^-$, and $\rho$ represent
relative changes in demand loads and power generation, so that
constraints \eqref{eq:feas.rest.power.pv.p},
\eqref{eq:feas.rest.power.pq.p}, and \eqref{eq:feas.rest.power.pq.q}
represent power flow equations \eqref{eq:power.flow} in which the
loads $P_\cG$, $P_\cD$, and $Q_\cD$ are modified.  The parameters
$\omega_i$ represent positive weights on the changes to loads and
generation, indicating the desirability or undesirability of changes
to that node. We note the following points.
\bi
\item The same variable $\rho_i$ is used in the active and reactive
  power balance equations \eqref{eq:feas.rest.power.pq.p} and
  \eqref{eq:feas.rest.power.pq.q}, since active and reactive load
  shedding should occur in the same fraction.
\item Bound constraints \eqref{eq:feas.rest.power.sigma.p},
  \eqref{eq:feas.rest.power.sigma.m}, and
  \eqref{eq:feas.rest.power.rho} on the load shedding variables limit
  the adjustments to a reasonable range (which may be zero for some
  buses).
\item The weights $\omega_i$ could be set to large positive values to
  discourage changes on that node, and to smaller values when change
  is acceptable. The case in which no change at all is allowable on
  that node can be handled by setting the upper bound to zero in
  \eqref{eq:feas.rest.power.sigma.p},
  \eqref{eq:feas.rest.power.sigma.m}, or
  \eqref{eq:feas.rest.power.rho}. Throughout the paper, we assume that
  $\omega_i =1$ for all $i$, but note  that other positive values
  of these weights can be used without any complication to the model.
\item Power generation at the generator nodes may be either increased
  or decreased in general, but the loads at demand nodes can only
  decrease. (Upper bounds $\ub{\sigma}_i^+$, $\ub{\sigma}_i^-$, and
  $\ub{\rho}_i$ should not exceed 1. This means that the type of a bus
  --- generator or demand bus --- cannot be changed.)
\item The bounds \eqref{eq:feas.rest.power.v} guarantee that voltage
  levels are operationally viable.
\ei
The objective to be minimized in \eqref{eq:feas.rest.power} is the
weighted sum of power adjustments that are necessary to restore
feasibility to the power flow equations.
We define $\cF_L(\gamma)=+\infty$ when it is not possible to restore
feasibility by adjusting loads and generations (which usually happens
because the constraints regarding acceptable voltage levels
\eqref{eq:feas.rest.power.v} cannot be satisfied even when load
shedding is allowed).

The feasibility restoration problem \eqref{eq:feas.rest.power} is a
nonconvex smooth constrained optimization problem in general, so we
can expect to find only a local solution when using standard
algorithms for such problems.  The problem generalizes
\eqref{eq:power.flow} in that if a solution of the latter problem
exists, it will yield a global solution of \eqref{eq:feas.rest.power}
with an objective of zero when we set $\sigma^+_i=\sigma^-_i=0$ for $i
\in \cG$ and $\rho_i=0$ for $i \in \cD$, provided the voltage
constraints \eqref{eq:feas.rest.power.v} are satisfied.  Moreover, by
the well-known sparsity property induced by $\ell_1$ objectives, we
expect few of the components of $\sigma_{\cG}^+$, $\sigma_{\cG}^-$,
and $\rho_\cD$ to be nonzero at a typical solution of
\eqref{eq:feas.rest.power}. The problem \eqref{eq:feas.rest.power} may
also have operational relevance, guiding the grid operator toward a
set of decisions that can restore stable operation of the grid with
minimal disruption.

For convenience of later discussion, we state
\eqref{eq:feas.rest.power} in the following more compact form:
\bseq \label{eq:feas.rest}
\begin{align}
\cF_L(\gamma) :=
\omin{x,y} & p^Ty\\
\ost       & F_L(x,y;\gamma) = 0\\
           & \lb{x} \le x \le \ub{x}\\
           & 0 \le y \le \ub{y},
\end{align}
\eseq
where $F_L(x,y;\gamma)=0$ represents the equality constraints
\eqref{eq:feas.rest.power.pv.p}-\eqref{eq:feas.rest.power.pq.q}, $x$
includes the circuit variables $V$ and $\theta$, and $y$ includes the
power-adjustment variables $\sigma^+$, $\sigma^-$, and $\rho$.


\subsubsection{Bilevel Formulation}\label{sec:problem.vli}

The bilevel optimization formulation seeks the attack $\gamma$ for
which the power-adjustment objective $\cF_L$ is maximized subject to the
same attack budget constraints as in \eqref{eq:model.volt}, that is,
\bseq \label{eq:model.LS}
\begin{alignat}{1}
  \cH_L(\kappa,\ub\gamma) :=
  \omax{\gamma} & \cF_L(\gamma) \\
  \ost & \mathrlap{e^T\gamma \le \kappa\ub{\gamma}} \label{eq:model.LS.total}\\
       & \mathrlap{0 \le \gamma \le \ub{\gamma}e.}  \label{eq:model.LS.budget}
\end{alignat}\eseq
By substituting from \eqref{eq:feas.rest}, we obtain a max-min
problem:
\bseq \label{eq:model.LS.alter}
\begin{alignat}{1}
  \cH_L(\kappa,\ub{\gamma}) :=
  \max_{\gamma}\,\omin{x,y} & p^Ty \\
  \ost & F_L(x,y;\gamma) = 0 \\
       & \lb{x} \le x \le \ub{x} \\
       & 0 \le y \le \ub{y} \\
       & e^T \gamma \le \kappa\ub{\gamma} \\
       & 0 \le \gamma \le  \ub{\gamma}e.
\end{alignat}\eseq

Bilevel optimization problems are, in general, difficult to solve.
For problems of the form \eqref{eq:model.LS.alter}, it is possible for
the upper-level objective $\cF_L$ to change discontinuously at some
values of $\gamma$, even when the constraint function $F_L$ is smooth
and nonlinear.

For the power-adjustment formulation, there is an additional
complication: For most feasible values of $\gamma$, the objective
function is zero. This is because power grids are often robust to
small perturbations, so when even when many impedances change, it is
often possible to continue meeting all demands while respecting
operational limits on the voltage values.  This feature makes it
difficult to search for the optimal $\gamma$, since it is difficult
even to find a starting value of $\gamma$ that causes nonzero
disruption.  We have developed specialized heuristics to address this
issue; these are described in Section~\ref{sec:alg.LS.init}.

\section{Algorithm Description}\label{sec:alg}

We discuss a first-order method for the following formulation, which
generalizes \eqref{eq:model.volt} and \eqref{eq:model.LS}:
\bseq\label{eq:model.general}
\begin{align}
\label{eq:model.general.1}
\cH(\kappa,\ub\gamma) :=
  \omax{\gamma} & \cF(\gamma) \\
\label{eq:model.general.2}
  \ost & e^T \gamma \le \kappa\ub{\gamma}\\
\label{eq:model.general.3}
       & 0 \le \gamma \le \ub{\gamma}e.
\end{align}
\eseq
Although the objective $\cF$ is not convex or smooth, we solve it with
the classical Frank-Wolfe method (also known as the conditional
gradient method), which we describe in the next subsection.

\subsection{Frank-Wolfe Algorithm} \label{sec:alg.FW}

The Frank-Wolfe algorithm \cite{FraW56} solves a sequence of
subproblems in which a first-order approximation to the objective
around the current iterate is minimized over the given feasible set.
If the objective $\cF$ in \eqref{eq:model.general} were smooth, we
would solve the following problem at the $k$th iterate $\gamma^k$:
\beq\label{eq:frank.wolfe.lp}
\begin{aligned}
w^k := \, \oargmax{w} & (g^k)^T(w-\gamma^k) \\
          \ost        & e^Tw \le \kappa\ub{\gamma}, \;\;
                          0 \le w \le \ub{\gamma},
\end{aligned}\eeq
where $g^k$ is a gradient $\cF(\gamma)$ at $\gamma^k$. The new iterate is
obtained by setting
\[
  \gamma^{k+1} = \gamma^k + \alpha_k (w^k-\gamma^k),
\]
for some $\alpha_k \in (0,1]$. (Frank and Wolfe~\cite{FraW56} give a
  specific formula for $\alpha_k$ that guarantees a sublinear
  convergence rate for smooth convex $\cF$. An exact line search would
  yield a similar rate.)
Because of the special nature of our constraint set, the problem
\eqref{eq:frank.wolfe.lp} is a linear program with a closed-form
solution, whose components $w^k_i$ are defined as follows:
\beq  \label{eq:expr.wk}
w^k_i =
 \begin{cases}
   \ub{\gamma} & \mbox{if $g^k_i$ is one of $\kappa$ largest entries in } g^k,\\
   0 & \mbox{otherwise.}
 \end{cases}
\eeq
We determine the step size $\alpha_k$ by a standard backtracking
procedure. Given a constant $\phi \in (0,1)$, and starting from
$\alpha=1$, we decrease the step size by $\alpha \leftarrow \phi
\alpha$ until the following sufficient decrease condition is satisfied
for a small $c_1\in(0,1)$.
\begin{equation}
  \cF(\gamma^k+\alpha (w^k-\gamma^k)) \ge \cF(\gamma^k) + c_1\alpha{g^k}^T(w^k-\gamma^k).
\end{equation}
We define $\alpha_k$ to be the value of $\alpha$ accepted by this
criterion. The algorithmic framework is shown in
Algorithm~\ref{alg:VA}. We terminate when the step $(w^k-\gamma^k)$
becomes small, or when the step size $\alpha$ becomes less than a
predefined $\alpha_{\min}>0$.

Convergence behavior of the Frank-Wolfe procedure with backtracking
line search for the smooth nonconvex case has been analyzed by
Dunn~\cite[Theorem~4.1]{Dunn80a}, where it is shown that accumulation
points are stationary. (This result does not apply directly to our
cases, because of potential nonsmoothness of the objectives.)

\begin{algorithm}[!h]\footnotesize
\caption{$\mbox{\sc Vulnerability Analysis}$}\label{alg:VA}
\begin{algorithmic}[1]
\Require
  \Statex $\ub{\gamma}$: Upper bound for impedance increases $\gamma_i$, $i \in \cL$;
  \Statex $\kappa$: Number of lines to attack;
  \Statex $\gamma^0$: Feasible initial value of $\gamma$;
\Ensure
  \Statex $\gamma^*$: Impedance vector that optimizes the attack;
\medskip
\State $k\gets 0$;
\While {$k\le ${\sc MaxIter}}
  \State Find the gradient $g^k$ of the objective $\cF$ at $\gamma^k$;
  \State Find linearized optimum $w^k$ from \eqref{eq:expr.wk};
  \State Use the backtracking to find step size $\alpha_k \in [0,1]$;
  \State {$\gamma^{k+1}\gets \gamma^k+\alpha_k (w^k-\gamma^k)$};
  \State $k\gets k+1$;
  \State Stop if termination conditions are satisfied, and set $\gamma^*\gets \gamma^k$;
\EndWhile
\end{algorithmic}
\end{algorithm}

\subsection{Gradient Calculation} \label{sec:gradient}

Algorithm~\ref{alg:VA} requires calculation of a gradient $g^k$ of the
objective function $\cF$ at the current iterate $\gamma^k$. We have
noted already that the power-adjustment objective $\cF = \cF_L$ may be
nonsmooth, due to changes in the active set of the subproblem
\eqref{eq:feas.rest}, so the gradient may not be well defined.  We
note however that $\cF_L$ can reasonably be assumed to be smooth
almost everywhere; changes to the active set can be expected to happen
only on a set of measure zero in the feasible space for $\gamma$.  Our
algorithm does not appear to encounter values of $\gamma$ where
$\cF_L$ is nondifferentiable in practice.

We outline a scheme for calculating gradients of $\cF_V$ and $\cF_L$
in Appendix~\ref{app:gradient}. The technique is essentially to use
the implicit function theorem to find sensitivities of the variables
in the problems that define $\cF_V$ and $\cF_L$ to the parameters
$\gamma$, around the current solution of these problems, and then
proceed to find the sensitivities of the optimal objective value for
these problems to $\gamma$.

\subsection{Power-Adjustment Model Initialization}\label{sec:alg.LS.init}

\begin{figure}\footnotesize
\centering \subfloat[Even distribution of relative impedance changes]{%
  \includepgfplots[width=.95\linewidth, axisratio=4]{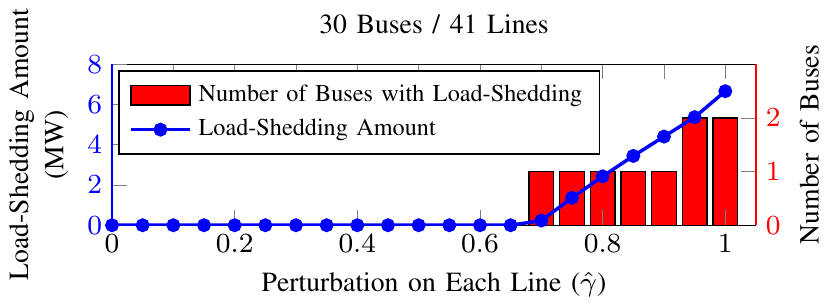}%
  \label{fig:30.perturb.evenly}}\\
\subfloat[``Safe'' distribution of relative impedance changes, to minimize total power adjustment]{%
  \includepgfplots[width=.95\linewidth]{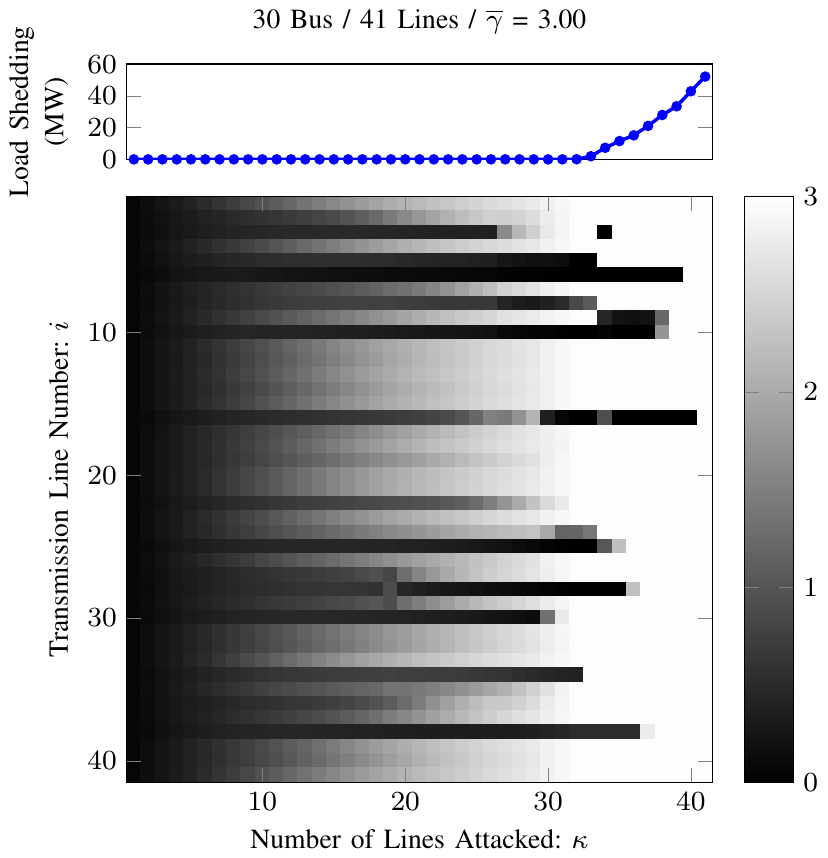}%
  \label{fig:30.perturb.safely}}
\caption{Power adjustment as a result of changes to line impedances on
  the 30-Bus data set. It is often the case that load shedding is not
  required even when the disturbance to the system is quite
  substantial.}\label{fig:30.perturb}
\end{figure}

As mentioned above, the objective value of the bilevel formulation is
zero for most feasible values of $\gamma$. It tends to be nonzero only
on parts of the feasible region defined by \eqref{eq:model.general.2},
\eqref{eq:model.general.3} that correspond to near-maximal attacks
focused on small numbers of buses.

To illustrate this point, we perform experiments on the 30-Bus case
({\tt case30.m} from \MATPOWER, originally from \cite{AlsS74}) in
which we monitor the power-adjustment objective $\cF_L(\gamma)$
in \eqref{eq:feas.rest.power} as the impedances are increased. In
Figure~\subref*{fig:30.perturb.evenly}, we plot $\cF_L(\gamma)$ in for
the values $\gamma = \hat{\gamma} e$, where $\hat{\gamma}$ is a
nonnegative scalar parameter that is increased progressively from $0$
to $1$. That is, impedances are increase {\em evenly across all
  transmission lines}. Note that $\cF_L(\hat{\gamma}e)$ is zero for
$\hat{\gamma} \in [0,0.65]$, while for $\hat{\gamma} > 0.65$, load
shedding occurs on one or two demand buses. This observation implies
that any value $\gamma$ along the line $\hat{\gamma} e$ (for
$\hat{\gamma} \in [0,.65]$) is a global minimizer of the bilevel
problem \eqref{eq:model.LS}. The gradient is zero at each of these
points, so optimization methods that construct the search direction
from gradients cannot make progress if started anywhere along this
line (or indeed from anywhere in a large neighborhood of this line).

If we are allowed to distribute a ``budget'' of impedance increases
unequally between lines, so as to minimize the total amount of
power adjustment required, even greater disturbances can be tolerated.
To describe this greater tolerance, we consider the following problem
that is sliglty modified from \eqref{eq:model.LS.alter}:
\bseq \label{eq:model.safe}
\begin{alignat}{2}
 \cH_S(\kappa,\ub{\gamma}) := \omin{x,y,\gamma} & p^Ty\\
  \ost & F_L(x,y;\gamma) = 0\\
       & \lb{x} \le x \le \ub{x}\\
       & 0 \le y \le \ub{y}\\
       & e^T \gamma = \kappa\ub{\gamma} \label{eq:model.safe.budget}\\
       & 0 \le \gamma \le \ub{\gamma}e. \label{eq:model.safe.ub}
\end{alignat}\eseq
Note that \eqref{eq:model.safe} is different from \eqref{eq:model.LS}
in two respects: (a) it is a single-level minimization problem whose
variables are $x, y,$ and $\gamma$; and (b) the budget is enforced
with the equality constraint \eqref{eq:model.safe.budget}. Thus this
problem finds a ``safe'' way to distribute the fixed budget
($\kappa\ub{\gamma}$) to transmission lines while the total
load-shedding $p^Ty$ is minimized.  We solved this problem for an
upper bound $\ub{\gamma}=3$ with $\kappa$, which is increased
progressively from $1$ to $41$. The top chart in
Figure~\subref*{fig:30.perturb.safely} shows that it is possible to
increase $\kappa$ to about $33$ before any load shedding takes place
at all. The lower chart depicts how the impedance changes are
distributed around the 41 lines in the grid, at the solution of
\eqref{eq:model.safe}, for each value of $\kappa$. Darker bars on the
graph show lines that can tolerate only a relatively small increase in
impedance before causing load shedding somewhere in the grid. The
lighter bars are those that can tolerate their impedance value
$\gamma_i$ being set to a value at or near the upper bound $3$ without
affecting load shedding. As an example: When $\kappa = 40$, we have
$\cH_S(40,3) \approx 55$, and the $\gamma$ that achieves this
power-adjustment value has components of $3$ on all lines except line
$16$, where it is zero.

The methodology used to derive Figure~\subref*{fig:30.perturb.safely}
can be used as a heuristic to identify a set of ``safe'' lines $\cS$
(whose impedances can be increased without affecting grid performance)
and a complementary set of ``vulnerable'' lines $\cW$ (for which
impedance increases are likely to lead to load shedding). In
Appendix~\ref{app:safe.vulnerable}, we describe the ESL (``estimating
safe lines'') procedure, Algorithm~\ref{alg:ESL}, for determining the
sets $\cS$ and $\cW$. Once we have determined the vulnerable lines
$\cW$, we define an impedance perturbation vector $\gamma'$ with the
following components:
\begin{equation} \label{eq:def.gammap}
\gamma'_i = \begin{cases} \ub{\gamma} & \;\; i \in \cW \\
0 & \;\; i \notin \cW,
\end{cases}
\end{equation}
where $\ub{\gamma}$ is the given upper bound on impedance on a given
line.  We then evaluate $\cF_L(\gamma')$ from
\eqref{eq:feas.rest.power}. If a node does not require any load
shedding under this maximal-perturbation setting, it is unlikely that
any attack on the vulnerable lines will lead to load shedding on this
node. We gather the other nodes --- those for which $\rho_i>0$ at the
solution of \eqref{eq:feas.rest.power.rho} with $\gamma=\gamma'$ ---
into a set $\cT$, the ``target nodes.''  Further explanation of the
definition of $\cT$ is given in Appendix~\ref{sec:alg.node}.

\begin{figure}\centering\footnotesize
 \subfloat[Target node]{%
  \includepgfplots[width=.8\linewidth]{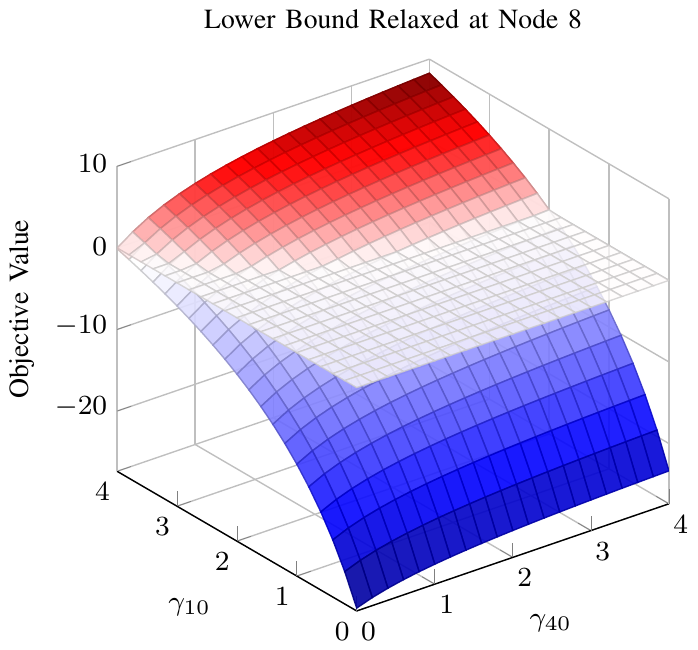}%
  \label{fig:target.30.T}}\ \\
 \subfloat[Non-target node]{%
  \includepgfplots[width=.8\linewidth]{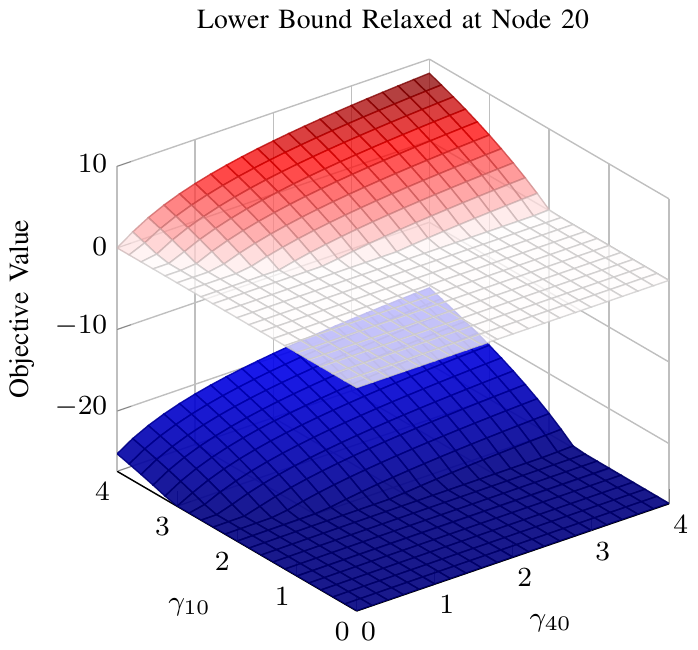}%
  \label{fig:target.30.NT}}
\caption{Extending the range of $\cF_L(\gamma)$ by allowing additional
  load on demand nodes. Adding load to target nodes (top figure)
  produces a useful extension of the range of the objective, so that
  its gradient yields a promising search direction for the
  maximization problem. Adding load to non-target nodes (bottom
  figure) simply shifts the objective down by a constant, so that the
  gradient in the flat region still does not yield a useful search
  direction for the maximization problem.}\label{fig:target.30}
\end{figure}

We use the target nodes to define a modification of the objective
$\cF_L$ that has the effect of shifting the {\em range} of the
function, in a way that makes gradient information relevant even at
values of $\gamma$ for which no power adjustments are required. The
idea is illustrated in Figure~\ref{fig:target.30}, where we show the
power-adjustment requirement on two different nodes (nodes 8 and 20)
of the 30-Bus system as a function of the values of two impedance
parameters --- those corresponding to lines 10 and 40.  In both graphs
of Figure~\ref{fig:target.30}, the top surfaces (shaded white and red)
represent the objective $\cF_L$ as a function of various values of the
pair $(\gamma_{10},\gamma_{40})$.  Note that $\cF_L$ takes the value
zero over much of the domain, but becomes positive when both
$\gamma_{10}$ and $\gamma_{40}$ are high. The lower surfaces in each
graph show how $\cF_L$ changes when we modify the subproblem in
\eqref{eq:feas.rest.power} by {\em removing the zero lower lower bound
  on the load $\rho_i$} in \eqref{eq:feas.rest.power.rho}, where $i=8$
(a target node) in Figure~\subref*{fig:target.30.T} and $i=20$ (a
non-target node) in Figure~\subref*{fig:target.30.NT}. Removal of the
lower bound has the effect of allowing load to be {\em added} to the
node in question. This is not an action that would be
operationally desirable, but as we see from the blue surface in
Figure~\subref*{fig:target.30.T}, it changes the nature of $\cF_L$ in useful
ways. The effect of removing the lower bound on $\rho_8$
(Figure~\subref*{fig:target.30.T}) is to extend the range of $\cF_L$
so that its derivative at {\em any} point in the domain gives useful
information about a good search direction. In a sense, the
extended-range version appears to be a natural extension of the
original objective $\cF_L$. By contrast, removal of the lower bound on
$\rho_{20}$ (Figure~\subref*{fig:target.30.NT}) causes the function to
simply be shifted downward by a roughly constant amount for all pairs
of impedance perturbation values. This is because, being a non-target
node, increased load on this node can be met, even after the grid is
damaged by the impedance attack. We conclude that {\em removing lower
  bounds on $\rho_i$ for target nodes} $i \in \cT$ provides a
potentially useful extension of the range of the function $\cF_L$,
whereas the same cannot be said for non-target nodes.

Motivated by these observations, we modify Algorithm~\ref{alg:VA} as
follows. We start by removing all lower bounds in
\eqref{eq:feas.rest.power.rho} on target nodes $i \in \cT$.  At each
iteration of the algorithm, after taking a step, we check to see if
any of the $\rho_i$ obtained by solving the subproblem
\eqref{eq:feas.rest.power} at the latest iteration are negative. If
so, we {\em reset the lower bound on the most negative value of
  $\rho_i$ to zero}, before moving on to the next iteration. The
algorithm does not terminate until all $\rho_i$ are nonnegative. The
modified procedure is shown as Algorithm~\ref{alg:VA.LS}.

\begin{algorithm}[t]\footnotesize
\caption{{\sc Vulnerability Analysis: Power-Adjustment Model}}\label{alg:VA.LS}
\begin{algorithmic}[1]
\Require
  \Statex $\ub{\gamma}$: Upper bound for impedance increases $\gamma_i$, $i \in \cL$;
\Statex  $\kappa$: Number of lines to attack;
\Statex $\gamma^0$: Feasible initial value of $\gamma$;
\Ensure
  \Statex $\gamma^*$: Impedance vector that optimizes the attack;
\medskip
\State Find a set of vulnerable lines $\cW$ and target nodes $\cT$ using the ESL procedure (Algorithm~\ref{alg:ESL});
\State Set lower bound of $\rho_i$ (for target nodes $i\in\cT$) to
$-\infty$;
\State {$k\gets 0$;}
\While {$k\le ${\sc MaxIter}}
  \State Find  gradient $g^k$ of $\cF_L$ at $\gamma^k$.
  \State Find linearized optimum $w^k$ from \eqref{eq:expr.wk};
  \State Use the backtracking to find step size $\alpha_k \in [0,1]$;
  \State {$\gamma^{k+1}\gets \gamma^k+\alpha_k (w^k-\gamma^k)$};
  \State Identify $i$ such that $\rho_i = \arg \min_j \rho_j$, where $\rho_j$ are the
\Statex $\quad\quad$ power-adjustment variables from \eqref{eq:feas.rest.power};
  \If{$\rho_i<0$}
  \State reset lower bound on $\rho_i$ to zero;
  \EndIf
  \State $k\gets k+1$;
  \State Stop if terminating conditions (including $\rho_j \ge 0$ for all power-adjustment variables $\rho_j$) are satisfied;
\EndWhile
\State $\gamma^*\gets \gamma^k$;
\end{algorithmic}
\end{algorithm}

\section{Experimental Results}\label{sec:result}

We present the results obtained with our formulations and algorithms
on the IEEE 118-Bus system and Polish 2383-Bus system. Our
implementations use \MATLAB{}\footnote{Version 8.1.0.604 (R2013a)}
with \IPOPT{}\footnote{Version 3.11.7}~(W\"achter and
Biegler~\cite{WacB06}) as the nonlinear solver for evaluating $\cF_L$
\eqref{eq:feas.rest.power} in the power-adjustment (bilevel)
model. For the test case data and calculation of the electric circuit
parameters, the codes from \MATPOWER{}~\cite{ZimMT11} are used
extensively. The codes were executed on a Macbook Pro (2 GHz Intel
Core i7 processor) with 8GB RAM.

\begin{table}[!h]\centering\footnotesize
\setlength{\tabcolsep}{.35em}
\caption{Test Cases for Experiments}\label{tbl:testcase}
\begin{tabular}{|c||C{.7cm}|C{.7cm}|C{.7cm}|C{.7cm}||}\cline{2-5}
\mc{1}{c||}{} & \mc{4}{c||}{Test Cases} \\\cline{2-5}
\mc{1}{c||}{} & 1 & 2 & 3 & 4 \\\hline
Filename (in \MATPOWER{}) & \mc{2}{c|}{case118.m} & \mc{2}{c||}{case2383wp.m} \\\hline
Number of Nodes     & \mc{2}{c|}{118}       & \mc{2}{c||}{2383}  \\\hline
Number of Lines     & \mc{2}{c|}{186}       & \mc{2}{c||}{2896}  \\\hline
Number of Lines to Attack $\kappa$   & 3 & 5 & 3 & 5 \\\hline
Perturbation Upper Bound $\ub{\gamma}$ & \mc{2}{c|}{3}         & \mc{2}{c||}{2}      \\\hline
\tcell{Backtracking Parameters $(c_0, c_1, \alpha_{\min})$}  & \mc{4}{c||}{(0.5, 0.01, 0.01)} \\\hline
Voltage limits $(\lb{V}, \ub{V})$
       & \mc{2}{c|}{$(0.93, 1.07)$} & \mc{2}{c||}{$(0.89, 1.12)$} \\\hline
Line Screening Threshold $\eta$ & \mc{4}{c||}{0.9} \\\hline
\end{tabular}
\end{table}

Information about the test case instances and algorithmic parameters
are given in Table~\ref{tbl:testcase}. There are two instances for
each of the two grids, corresponding to 3-line and 5-line attacks,
respectively. The table shows voltage magnitude limits that are
applied in the power-adjustment model, together with the value of
$\eta$ that is used in the ESL procedure (Algorithm~\ref{alg:ESL} from
Appendix~\ref{sec:alg.line}).

In the power-adjustment model \eqref{eq:feas.rest.power}, the upper
bounds $\overline{\sigma}_i^+$, $\overline{\sigma}_i^-$, and
$\overline{\rho}_i$ on the power-adjustment variables are set to 1 for
most buses, thus allowing full load shedding. If a bus violates our
assumption on power injection --- that is, if $P_i\le0$ for bus
$i\in\cG$ or $P_i\ge0$ or $Q_i\ge0$ for bus $i\in\cD$ --- the
load-shedding upper bound for that bus is set to 0, disallowing power
adjustment on that bus. The power-adjustment objective $\cF_L$ is
considered to be nonzero if it is at least $10^{-3}$ megawatt (MW).

\subsection{Voltage Disturbance Model} \label{sec:vd}

We discuss first results obtained with the voltage disturbance model
\eqref{eq:def.FV}-\eqref{eq:model.volt} applied to the four test cases
of Table~\ref{tbl:testcase}.

\subsubsection*{118-Bus System}

\begin{figure}[h!]
\centering\footnotesize
\subfloat[Distribution of Voltage Magnitudes Before and After Attack.]{%
  \includepgfplots[width=.98\linewidth,axisratio=3]{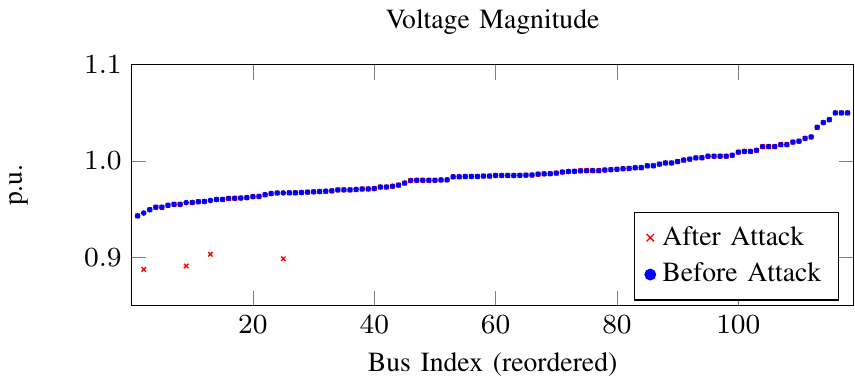}%
\label{fig:result.volt.118.k3.mag}
}\\
\subfloat[Distribution of Voltage Magnitudes Changes After Attack.]{%
  \includepgfplots[width=.98\linewidth,axisratio=3]{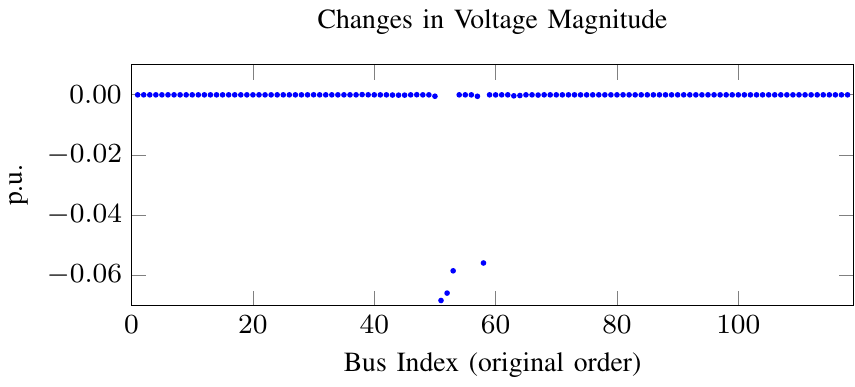}%
\label{fig:result.volt.118.k3.change}}
\caption{Voltage Disturbance Model: 118-Bus System with $\kappa=3$}\label{fig:result.volt.118.k3}
\end{figure}

\begin{figure*}\centering
\includegraphics[width=.65\linewidth]{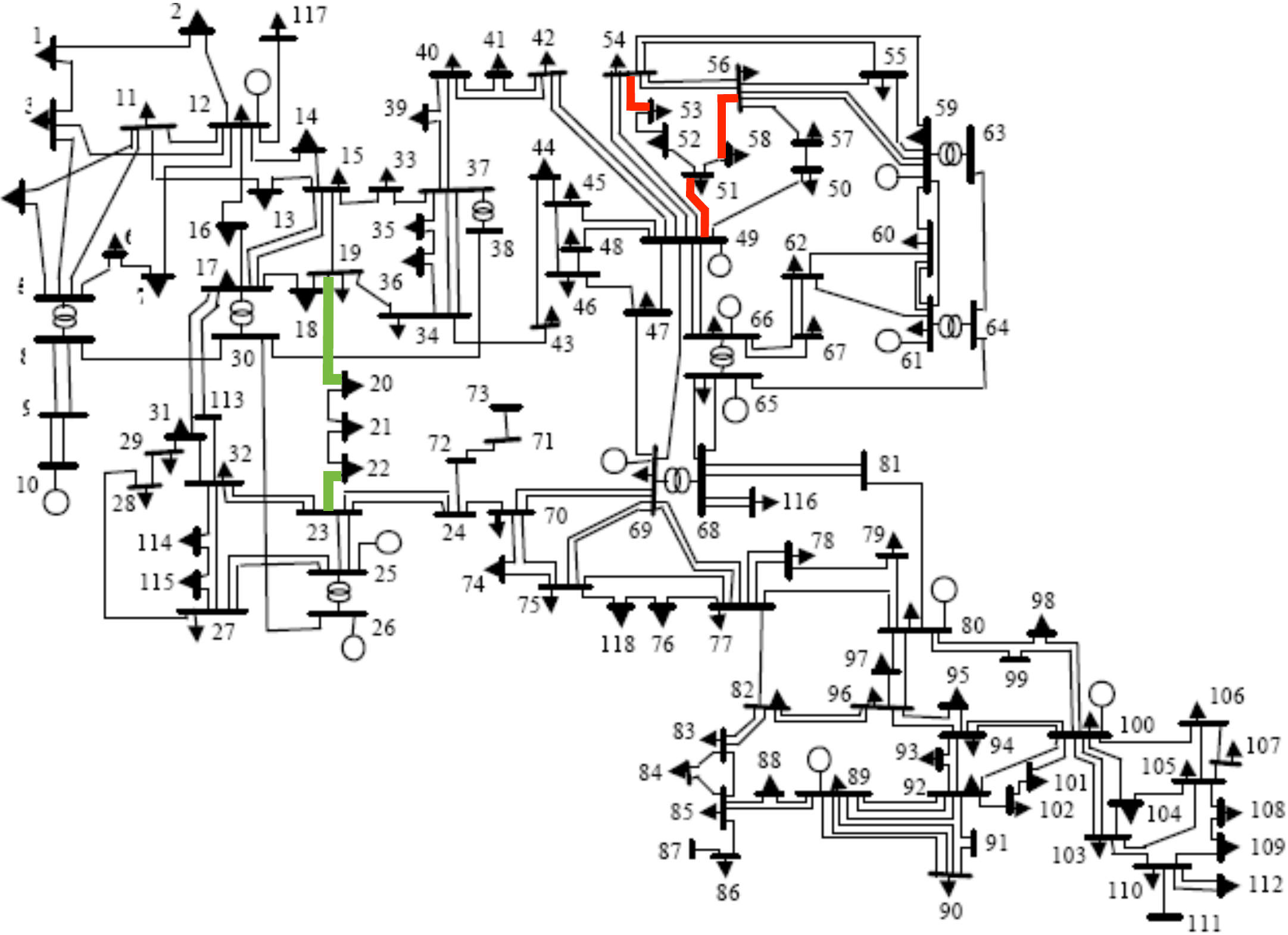}\vspace{-3mm}
\caption{Voltage Disturbance Model: Attacks on the 118-Bus System.  The
  transmission lines in red are are the optimal attack for
  $\kappa=3$. The lines in green are added to the optimal 3-line
  attack when $\kappa$ is increased to 5.}\vspace{-2mm}
\label{fig:result.volt.118}
\end{figure*}

For a 3-line attack problem on IEEE 118-Bus system ($\kappa=3$),
Algorithm~\ref{alg:VA} converges in 5 iterations and identifies
exactly three lines to attack with maximal impedance increase: lines
71, 74, and 82 (as shown in Table~\ref{tbl:result.volt.118.k3}).
As a result of this attack, voltage magnitudes at four buses decrease
significantly, by up to 0.07~p.u., as shown in
Figure~\ref{fig:result.volt.118.k3}. (In
Figure~\subref*{fig:result.volt.118.k3.mag}, the buses are reordered
in increasing order of voltage magnitude on the undisturbed system. In
Figure~\subref*{fig:result.volt.118.k3.change}, the buses are indexed
in their original order.) The attack is visualized in
Figure~\ref{fig:result.volt.118}, where we see that its effect is
essentially to isolate buses 51, 52, 53, and 58; the attacked lines
are colored in red.

\begin{table}[ht!]
\centering\footnotesize
\caption{Voltage Disturbance Model: 118-Bus System with $\kappa=3$}\label{tbl:result.volt.118.k3}
\captionsetup[subfloat]{labelformat=empty}\vspace{-5px}
\subfloat[Optimal Attack (as determined by our algorithm)]{%
\makebox[.95\linewidth]{%
\begin{tabular}{||c|c|c||c||}
\hline
Line & \mc{2}{c||}{Buses}  & Continuous \\
\cline{2-3}
No. & From &  To & Attack ($\gamma_i$) \\\hline\hline
 71 &  49  &  51  & 3.00  \\
 74 &  53  &  54  & 3.00  \\
 82 &  56  &  58  & 3.00  \\\hline\hline
\mc{3}{||c||}{Objective} & $4.04\times10^{-2}$ \\\hline
\end{tabular}}}
\end{table}
\begin{table}[t]
\centering\footnotesize
\caption{Voltage Disturbance Model: 118-Bus System with $\kappa=5$}
\label{tbl:result.volt.118.k5}
\captionsetup[subfloat]{labelformat=empty}\vspace{-5px}
\subfloat[Optimal Attack (as determined by our Algorithm)]{%
\makebox[0.95\linewidth]{%
\begin{tabular}{||c|c|c||c||}
\hline
Line & \mc{2}{c||}{Buses}  & Continuous \\
\cline{2-3}
No. & From &  To & Attack ($\gamma_i$) \\\hline\hline
 25 &  19  &  20 & 3.00 \\
 29 &  22  &  23 & 3.00 \\
 71 &  49  &  51 & 3.00 \\
 74 &  53  &  54 & 3.00 \\
 82 &  56  &  58 & 3.00 \\\hline\hline
\mc{3}{||c||}{Objective} & $5.03\times10^{-2}$ \\\hline
\end{tabular}}}
\end{table}
\begin{figure}[ht!]
\centering\footnotesize
\subfloat[Distribution of Voltage Magnitude Before and After Attack.]{%
  \includepgfplots[width=.98\linewidth,axisratio=3]{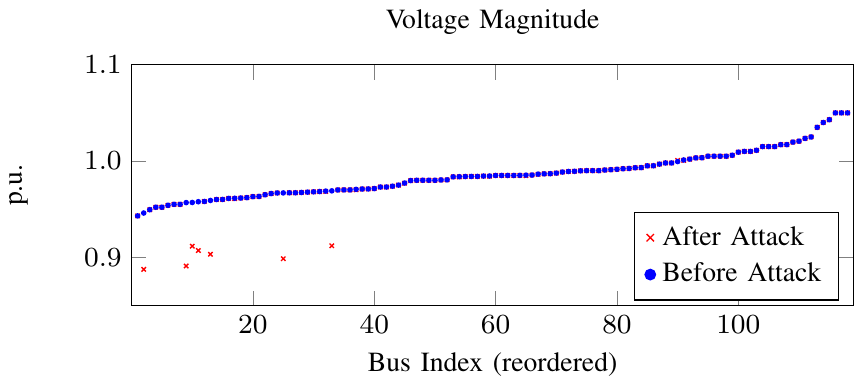}%
}\\
\subfloat[Distribution of Voltage Magnitude Changes and After Attack.]{%
  \includepgfplots[width=.98\linewidth,axisratio=3]{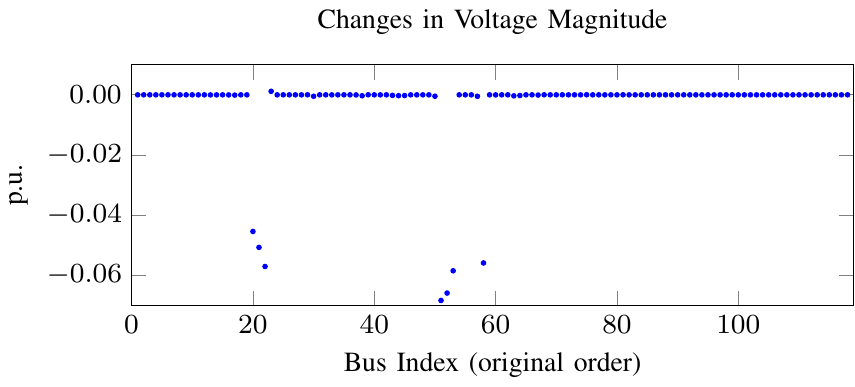}%
\label{fig:result.volt.118.k5.change}}
\caption{Voltage Disturbance Model: 118-Bus System with $\kappa=5$}
\end{figure}


For the second test instance, on the IEEE 118-Bus system with
$\kappa=5$, the algorithm identifies exactly five lines to attack at
the maximal impedance increase --- lines 25, 29, 71, 74, and 82 (see
Table~\ref{tbl:result.volt.118.k5}) --- which includes the three
lines identified in the first attack. With this stronger attack, there
is significant voltage drop on seven buses. Algorithm~\ref{alg:VA}
takes 8 iterations to converge to the solution.  As
Figure~\ref{fig:result.volt.118} shows, attacking the additional two
lines (colored in green) has the effect of creating another
``island,'' consisting of buses 20, 21 and 22. The additional voltage
drops seen Figure~\subref*{fig:result.volt.118.k5.change} are from these
buses.

\subsubsection*{2383-Bus System}

\begin{table}[h]
\centering\footnotesize
\caption{Voltage Disturbance Model: 2383-Bus System with $\kappa=3$}
\label{tbl:result.volt.2383.k3}
\captionsetup[subfloat]{labelformat=empty}\vspace{-10px}
\subfloat[Optimal Attack (as determined by our Algorithm)]{%
\begin{tabular}{||c|c|c||c|c|c||}
\hline
Line & \mc{2}{c||}{Buses}  & Continuous & Top-3 & Best-3 \\
\cline{2-3}
No. & From &  To & Attack ($\gamma_i$) & Attack & Attack \\\hline\hline
  5 &  10  &   3 & 1.04 & 2.00 & 2.00 \\
404 & 434  & 188 & 0.25 &      &      \\
405 & 437  & 188 & 2.00 & 2.00 & 2.00 \\
467 & 340  & 218 & 2.00 & 2.00 & 2.00 \\
501 & 340  & 240 & 0.71 &      &      \\\hline\hline
\mc{3}{||c||}{Objective} & 0.514 & 0.501 & 0.501 \\\hline
\end{tabular}}
\end{table}
\begin{figure}[h!]
\centering\footnotesize\vspace{0mm}
\subfloat[Distribution of Voltage Magnitude Before and After Attack. (Best-3)]{%
  \includepgfplots[width=.98\linewidth,axisratio=3]{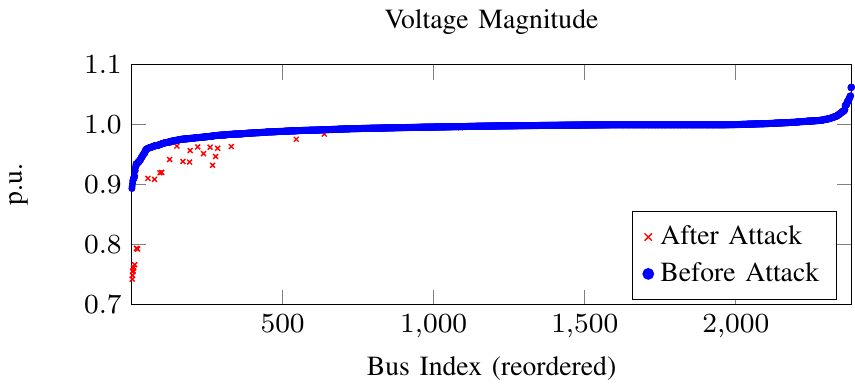}}\\
\subfloat[Distribution of Voltage Magnitude Changes After Attack. (Best-3)]{%
  \includepgfplots[width=.98\linewidth,axisratio=3]{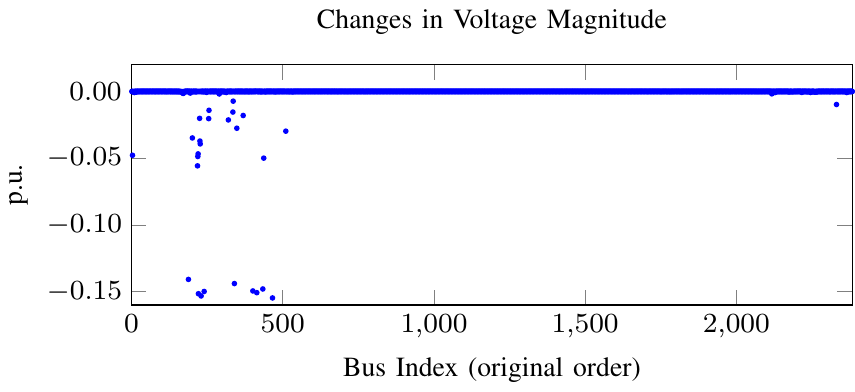}%
\label{fig:result.volt.2383.k3.change}}
\caption{Voltage Disturbance Model: 2383-Bus System with $\kappa=3$}\label{fig:result.volt.2383.k3}\vspace{-8mm}
\end{figure}

Results for our third test instance in Table~\ref{tbl:testcase}, which
considers the 2383-Bus model with attack limit defined by $\kappa=3$,
are shown in Table~\ref{tbl:result.volt.2383.k3}. The continuous
impedance attack is distributed into 5 lines, identified after ten
iterations of Algorithm~\ref{alg:VA}. (The lines involved in the
attack are revealed at iteration five, while the remaining five
iterations make minor adjustments to the impedance values.)

This 5-line solution can be used to identify the most disruptive set
of three lines using two heuristics: (a) choose the three lines $i$
such that $\gamma_i$ are one of the largest three entries in $\gamma$
(called the \textit{Top-3 attack}); and (b) try all possible 3-line
combinations of the five lines with highest impedance (the
\textit{Best-3 attack}). For this specific case, the Top-3 attack and
Best-3 attack are the same, consisting of lines 5, 405, and 467. These
3-line attacks, with impedance set to their maximum values on all
three lines, gives a slightly smaller objective value than the
continuous attack.

Figure~\subref*{fig:result.volt.2383.k3.change} shows how the
voltage magnitude changes when the Best-3 and Top-3 attacks are
executed.  Similarly to 118-Bus cases, there is a relatively small
number of buses in which the voltage drops significantly from its
original value, but large voltage changes are seen on some buses,
with some voltage magnitudes below 0.8~p.u.

For our fourth test instance in Table~\ref{tbl:testcase}, for which
$\kappa=5$ on the 2383-Bus system, two iterations of
Algorithm~\ref{alg:VA} suffice to identify an attack (on lines 404,
405, 467, 479, and 501) that makes the power flow problem infeasible,
that is, there is no $(V,\theta)$ that satisfies
$F(V,\theta;\gamma)=0$ under this attack. Hence, unless the grid
operator takes action (to change loads or generator outputs, for
example), an attack on these five lines renders the grid inoperable.

\subsubsection*{Comparison with Continuous Optimization Model}

In Section~\ref{sec:problem.volt}, we mentioned that the voltage
disturbance model can be written as a continuous optimization model
\eqref{eq:model.volt.nlp}, except that the latter model does not
handle infeasibility appropriately. We verify the properties of this
alternative formulation by solving it with the nonlinear
interior-point solver \IPOPT{}.  We find that in the first three
instances of Table~\ref{tbl:testcase}, the solutions obtained from
\eqref{eq:model.volt} match those we described above. In the fourth
test case --- the 2383-Bus model with $\kappa=5$ --- the model
\eqref{eq:model.volt.nlp} identifies a solution that maximizes
disruption {\em subject to the power flow equations
  \eqref{eq:power.flow} remaining feasible}. Our model
\eqref{eq:model.volt}, which detects infeasibility of the grid under
an attack of this strength, yields the more informative outcome.

\subsection{Power-Adjustment Model}

We now present results for the power-adjustment model
\eqref{eq:feas.rest.power}, \eqref{eq:model.LS}, for the four test
instances of Table~\ref{tbl:testcase}.

\subsubsection*{118-Bus System}
The ESL procedure (Algorithm~\ref{alg:ESL}), which is described in
Section~\ref{sec:alg.LS.init} and Appendix~\ref{sec:alg.line}, is
applied to the 118-Bus system to identify vulnerable lines $\cW=\{29,
71, 74, 96, 184\}$ and target nodes $\cT=\{52, 53, 117\}$. We use
these settings of $\cW$ and $\cT$ in Algorithm~\ref{alg:VA.LS} to
solve the first two instances in Table~\ref{tbl:testcase}.

For the first test case ($\kappa=3$), Algorithm~\ref{alg:VA.LS}
terminates after eight iterations, at the attack shown in
Table~\ref{tbl:result.load.118.k3}. Four lines are involved in this
attack; we reduce to three-line attacks ``Top-3'' and ``Best-3'' as in
 Subsection~\ref{sec:vd}. The Top-3 and Best-3 attacks coincide (and
are the same as those obtained for the voltage disturbance model in
Subsection~\ref{sec:vd}) and have a slightly smaller objective value
than the continuous attack. We note that in these 3-line attacks, the
voltage magnitudes of buses 52 and 53 are at their lower bound
$\lb{V}=.93$, and load shedding is required for buses 51 and 53, the
total amount of load shedding being 22.13 MW.

\begin{table}
\centering\footnotesize
\caption{Power-Adjustment Model: 118-Bus System with $\kappa=3$}
\label{tbl:result.load.118.k3}
\begin{tabular}{||c|c|c||c|c|c||}
\hline
Line & \mc{2}{c||}{Buses}  & Continuous & Top-3 & Best-3 \\
\cline{2-3}
No. & From &  To & Attack ($\gamma_i$) & Attack & Attack \\\hline\hline
 71 &  49  &  51 & 3.00 & 3.00 & 3.00 \\
 72 &  51  &  52 & 0.26 &      &      \\
 74 &  53  &  54 & 3.00 & 3.00 & 3.00 \\
 82 &  56  &  58 & 2.74 & 3.00 & 3.00 \\\hline\hline
\mc{3}{||c||}{Objective (MW)}            & 22.39      & 22.13  & 22.13 \\\hline
\mc{3}{||c||}{\tcell{Buses with \\ Load Shedding}}   & 51, 52, 53 & 51, 53 & 51, 53 \\\hline
\mc{3}{||c||}{\tcell{Buses at \\ Voltage Boundary}} & 52         & 52, 53 & 52, 53 \\\hline
\end{tabular}
\end{table}
\begin{table}
\centering\footnotesize
\caption{Five Most Effective Attacks from Exhaustive Enumeration: 118-Bus System with $\kappa=3$}
\label{tbl:result.load.comb.118.k3}
\begin{tabular}{||ccc|c||}\hline
\mc{3}{||c|}{Lines Selected} &  Power Adjustment (MW)\\\hline\hline
 71 &  74 &  82 & 22.13 \\
 71 &  72 &  74 & 19.07 \\
 71 &  74 &  83 & 17.21 \\
 71 &  74 & 184 & 15.61 \\
 71 &  74 &  97 & 13.27 \\\hline
\end{tabular}
\end{table}

In Table~\ref{tbl:result.load.comb.118.k3}, we show the top five
three-line attacks obtained by enumerating all $266,916$ possible
three-line attacks on this grid. The most disruptive attack is indeed
the one found by our algorithm.

\begin{table}
\centering\footnotesize
\caption{Power-Adjustment Model: 118-Bus System with $\kappa=5$}
\label{tbl:result.load.118.k5}
\begin{tabular}{||c|c|c||c|c|c||}
\hline
Line & \mc{2}{c||}{Buses}  & Continuous & Top-5 & Best-5 \\
\cline{2-3}
No. & From &  To & Attack ($\gamma_i$) & Attack & Attack \\\hline\hline
 71 &  49  &  51 & 3.00 & 3.00 & 3.00  \\
 72 &  51  &  52 & 0.36 &      & 3.00  \\
 74 &  53  &  54 & 3.00 & 3.00 & 3.00  \\
 75 &  49  &  54 & 0.74 &      &       \\
 76 &  49  &  54 & 0.81 & 3.00 &       \\
 81 &  50  &  57 & 0.29 &      &       \\
 82 &  56  &  58 & 3.00 & 3.00 & 3.00  \\
 97 &  64  &  65 & 0.80 &      &       \\
184 &  12  & 117 & 3.00 & 3.00 & 3.00  \\\hline\hline
\mc{3}{||c||}{Objective (MW)}            & 27.16 & 25.79  & 26.87 \\\hline
\mc{3}{||c||}{\tcell{Buses with \\ Power Adjustment}}   & 51, 52, 53, 117 & 51, 53, 117 & 52, 53, 117 \\\hline
\mc{3}{||c||}{\tcell{Buses at \\ Voltage Boundary}} & 52, 117 & 52, 53, 117 & 52, 117 \\\hline
\end{tabular}
\end{table}

For the second test instance (with $\kappa=5$), Algorithm~\ref{alg:VA}
requires twelve iterations to converge,
and distributes the impedance increases around nine lines, with the
maximum perturbation ($\gamma_i=3$) on four of them; see
Table~\ref{tbl:result.load.118.k5}.
The Top-5 and Best-5 attacks each involve the four lines with maximum
impedance increases, but differ in their choice of additional
line. Three buses require load shedding in the Best-5 attack, compared
to two in the Best-3 attack.

Unlike the previous results for $\kappa=3$, which target the same
lines as the voltage disturbance model, the Best-5 attack identified
in the power-adjustment model is slightly different from the Best-5
attack for the voltage-disturbance model. The power-adjustment attack
targets the region around buses 51, 52, and 53, as before, but also
``islands'' bus 117 (see Figure~\ref{fig:result.volt.118}), rather
than the buses 20, 21, and 22 that are attacked by the
voltage-disturbance model.

\subsubsection*{2383-Bus System}

For the Polish 2383-Bus system, the ESL procedure
(Algorithm~\ref{alg:ESL}) identified a set $\cW$ of 20 vulnerable
lines and, for upper bound $\ub{\gamma}=2$ on the impedance increase,
a set $\cT$ of 55 target nodes. (Algorithm~\ref{alg:ESL} required
about 270 seconds of run time on this instance.)

\begin{table}
\centering\footnotesize
\setlength{\tabcolsep}{.4em}
\caption{Power-Adjustment Model: 2383-Bus System with $\kappa=3$}
\label{tbl:result.load.2383.k3}
\begin{tabular}{||c|c|c||c|c|c||c|c||}
\hline
Line & \mc{2}{c||}{Buses}  & \mc{3}{c||}{From Bilevel Formulation} & \mc{2}{c||}{From $N-1$} \\
\cline{2-8}
No. & From &  To & Cont. ($\gamma_i$) & Top-3 & Best-3 & Top-3 & Best-3 \\\hline\hline
  5 &   10 &   3 &      &      &      & 2.00 &      \\
264 &  140 & 117 & 0.02 &      &      &      &      \\
268 &  126 & 118 & 2.00 & 2.00 & 2.00 & 2.00 & 2.00 \\
289 &  135 & 125 & 1.98 & 2.00 & 2.00 &      & 2.00 \\
296 &  145 & 128 & 2.00 & 2.00 & 2.00 & 2.00 & 2.00 \\\hline\hline
\mc{3}{||c||}{Objective (MW)}            & 577.80 & 577.75 & 577.75 & 303.07 & 577.75 \\\hline
\mc{3}{||c||}{\tcell{\# of Buses with \\ Power Adjustment}}   & \mc{3}{c|}{49} & 27 & 49\\\hline
\mc{3}{||c||}{\tcell{Buses at \\ Voltage Boundary}} & \mc{3}{c|}{145, 146, 1905} & \tcell{145, 146,\\ 230, 1905} & \tcell{145, 146,\\ 1905}\\\hline
\end{tabular}
\end{table}

Results for the power-adjustment model on or our third test case from
Table~\ref{tbl:testcase}, for the 2383-Bus system with $\kappa=3$, are
shown in
Table~\ref{tbl:result.load.2383.k3}. Algorithm~\ref{alg:VA.LS} returns
a four-line attack. Since the attack on one of these four lines (line
264) is negligible, we find that the Top-3 and Best-3 solutions both
attack lines 268, 289 and 296, with an active-power load shedding
577.75 MW, which is negligibly smaller than the optimal four-line
attack.  There are 49 buses which need load shedding for the
three-line attack, and three buses have voltage magnitudes at their
lower limits of $\lb{V}=.89$ --- another sign of stress on the grid.
The solution obtained for this test case is quite different from the
one from the voltage disturbance model. The lines 5, 404, 405, 467, and
501 which are identified by the voltage disturbance model (cf.
Table~\ref{tbl:result.volt.2383.k3}) also cause some load shedding
when they are attacked, but the effect is not as serious as the attack
on lines 268, 289, and 296.

Since it is computationally intractable to look at all possible
three-line combinations in a 2383-Bus grid, we apply an ``$N-1$
enumeration'' heuristic to explore the most promising part of the
space of three-line attacks.  In this heuristic, each line $i$ is
individually perturbed by setting $\gamma_i=\ub{\gamma}$, and we note
which of these perturbations require load shedding.  On this data set,
sixteen lines were identified as causing load shedding. We define a
``Top-3$(N-1)$'' attack to comprise the three lines that individually
cause the most load shedding, and the ``Best-3$(N-1)$'' attack to be the
most disruptive three-line combination drawn from these sixteen lines.
The resulting attacks are displayed in
Table~\ref{tbl:result.load.2383.k3}, alongside the Top-3 attack and
Best-3 attack. We see that the Best-3$(N-1)$ attack is identical to
the Top-3 and Best-3 attacks, while the Top-3$(N-1)$ attack is
inferior.

The optimal attacks for the power-adjustment model in this third test
instance are quite different from those obtained from the voltage
disturbance model, as we see by comparing
Tables~\ref{tbl:result.volt.2383.k3} and
\ref{tbl:result.load.2383.k3}.  We found, however, that if the lower
bound on voltage magnitude is changed from $\lb{V}=.89$ to
$\lb{V}=.87$, the optimal attack for the power-adjustment model is almost
identical to the voltage disturbance model. In the relaxed problem,
the buses 145, 146, and 1905 no longer have their voltage magnitude at
the lower bound, while buses 401 and 414 move to the relaxed lower
bound. The latter two buses are among those that suffer significant
voltage drop in the optimal voltage-disturbance attack of
Table~\ref{tbl:result.volt.2383.k3}.

\begin{table}
\centering\footnotesize
\setlength{\tabcolsep}{.30em}
\caption{Power-Adjustment Model: 2383-Bus System with $\kappa=5$}
\label{tbl:result.load.2383.k5}
\begin{tabular}{||c|c|c||c|c|c||c|c||}
\hline
Line & \mc{2}{c||}{Buses}  & \mc{3}{c||}{From Bilevel Formulation} & \mc{2}{c||}{From $N-1$} \\
\cline{2-8}
No. & From &  To & Cont. ($\gamma_i$) & Top-5 & Best-5 & Top-5 & Best-5 \\\hline\hline
  5   &   10 &    3  &      &      &      & 2.00 &      \\
  268 &  126 &  118  & 2.00 & 2.00 & 2.00 & 2.00 & 2.00  \\
  269 &  142 &  118  & 2.00 & 2.00 & 2.00 &      &       \\
  289 &  135 &  125  & 2.00 & 2.00 & 2.00 & 2.00 & 2.00  \\
  296 &  145 &  128  & 2.00 & 2.00 &      & 2.00 & 2.00  \\
  317 &  142 &  135  & 1.54 & 2.00 & 2.00 &      &       \\
  405 &  437 &  188  &      &      &      & 2.00 & 2.00  \\
  467 &  340 &  218  &      &      &      &      & 2.00  \\
 2142 & 1693 & 1658  & 0.46 &      & 2.00 &      &       \\\hline\hline
\mc{3}{||c||}{Objective (MW)} & 1109.72 & 1086.67 & 1460.25 & 594.09 & 597.71 \\\hline
\mc{3}{||c||}{\tcell{\# of Buses with \\ Power Adjustment}}  & 77 & 78 & 71 & 53 & 53\\\hline
\mc{3}{||c||}{\tcell{Buses at \\ Voltage Boundary}} & \mc{2}{c|}{145, 146, 1905} & 1905 & \tcell{145, 146,\\ 230, 414,\\ 434, 1905} & \tcell{145, 146,\\ 401, 414,\\ 1905}\\\hline
\end{tabular}
\end{table}

Table~\ref{tbl:result.load.2383.k5} shows results for our fourth test
instance from Table~\ref{tbl:testcase}, an optimal attack on the
2383-Bus system for $\kappa=5$. The attack determined by our procedure
is distributed around 6 lines. The Top-5, Best-5, Top-5$(N-1)$, and
Best-5$(N-1)$ attacks are calculated as described above.  The Best-5
attack is significantly more disruptive than the optimal continuous
attack identified by  Algorithm~\ref{alg:VA.LS}, probably because
of nonconcavity in the objective $\cF_L$. We note however that our
algorithm is much more useful in screening for the most disruptive
collection of lines than is the standard ``$N-1$'' screening
methodology: The Top-5 and Best-5 attacks are much more damaging than
the Top-5$(N-1)$ and Best-5$(N-1)$ attacks.


\section{Conclusions}\label{sec:conclusion}

We have proposed an attack model for assessing the vulnerability of
power grids. The attack consists in increasing the impedance of
transmission lines, with the resulting disruption to the grid measured
in two ways. The first technique is to observe changes in voltage
magnitudes at the buses; greater changes from the nominal values
indicate greater disruption. The second technique is to measure the
weighted sum of adjustments to load and generation that are needed to
restore stable operation of the grid, with the voltage magnitudes
confined to certain prespecified ranges. The two criteria give rise to
optimization problems with different properties. Both are solved with
a combination of known algorithms (such as Frank-Wolfe) and heuristics
that determine promising regions of the solution space.

In our computational results, we also use our algorithm as a screening
procedure for determining which collections of lines are likely to
cause the most disruptive attacks. By enumerating combinations of
lines from among those identified by our algorithms, we identify more
disruptive attacks than those produced by alternative screening
methods, such as the well known ``$N-1$'' criterion.

\appendix

\subsection{Computing Gradients} \label{app:gradient}

We describe here calculation of gradients for the functions
$\cF(\gamma)$ defined in Section~\ref{sec:problem}, which quantify the
grid disruption arising from an attack modeled by the impedance vector
$\gamma$. Both model functions considered here ---
\eqref{eq:def.FV} and \eqref{eq:feas.rest.power} --- have the
following form:
\begin{subequations} \label{eq:genF}
\begin{align}
\cF(\gamma) := \omin{z} & f(z) & \\
               \ost     & c_i(z,\gamma) = 0, \quad i=1,2,\dotsc,m, \\
                        & h_j(z) \ge 0, \quad j=1,2,\dotsc,r,
\end{align}
\end{subequations}
where the functions $f$, $c_i$, $i=1,2,\dotsc,m$, and $h_j$,
$j=1,2,\dotsc,r$ are all smooth. (Note that in our models, the
dependence on the upper-level variables $\gamma$ arises only through
the equality constraints, but our discussion can be extended without
conceptual difficulty to the case in which the objective and the
inequality constraints also depend on $\gamma$.)

We outline a technique for calculating the gradient $\nabla
\cF(\gamma)$.
We assume that the minimizing $z$ for the (generally nonconvex)
problem \eqref{eq:genF} has been identified and that it is denoted by
$z(\gamma)$. Moreover, we assume that $z(\gamma)$ is a {\em
  nondegenerate} solution of the minimization problem in
\eqref{eq:genF}. By this we mean that the linear independence
constraint qualification holds at the minimizer, that a strict
complementarity condition holds, and that second-order sufficient
conditions are satisfied at $z(\gamma)$. We show that under these
conditions, we can use the implicit function theorem to define the
gradient $\nabla \cF(\gamma)$ uniquely. While strong, these conditions
are not impractical; they appear to hold for all values of $z(\gamma)$
encountered by our algorithm, and it seems plausible that they would
hold for ``almost all'' values of $\gamma$.  The question of existence
of $\nabla \cF(\gamma)$ becomes much more complicated when these
conditions are not satisfied. When strict complementarity does not
hold, for example, the set of active inequality constraints is on the
verge of changing, an event usually associated with a point of
nonsmoothness of $\cF(\gamma)$.

The Karush-Kuhn-Tucker (KKT) conditions for optimality of $z$ in the
problem \eqref{eq:genF} for a given $\gamma$ are that there exist
scalars $\lambda_i$, $i=1,2,\dotsc,m$ and $\mu_j$, $j=1,2,\dotsc,r$
such that
\begin{subequations} \label{eq:kkt}
\begin{align}
\nabla f(z) - \sum_{i=1}^m \lambda_i \nabla_z c_i(z,\gamma) - \sum_{j=1}^r \mu_j \nabla h_j(z) = 0, \\
c_i(z,\gamma) = 0, \;\; i=1,2,\dotsc,m, \\
\mu_j  \ge 0, \;\; j \in \cA(z), \\
h_j(z) = 0, \;\; j \in \cA(z), \\
\mu_j = 0, \;\; j \in \{1,2,\dotsc,r \} \setminus \cA(z), \\
h_j(z) \ge 0, \;\; j \in \{1,2,\dotsc,r \} \setminus \cA(z),
\end{align}
\end{subequations}
where  the active set $\cA(z)$ is defined as follows:
\begin{equation} \label{eq:defA}
\cA(z) := \{j=1,2,\dotsc,r  \, : \, h_j(z)=0 \}.
\end{equation}
We use the following vector notation:
\[\begin{aligned}
c(z,\gamma) &= [c_i(z,\gamma)]_{i=1}^m, & \quad
\lambda &= [\lambda_i]_{i=1}^m,\\
\mu_\cA &= [\mu_j]_{j \in \cA}, & \quad
h_\cA(z) &= [h_j(z)]_{j \in \cA}.
\end{aligned}\]
The {\em linear independence constraint qualification (LICQ)} is:
\begin{equation}\label{eq:licq}
\parbox{.87\linewidth}{\centering
$\{ \nabla_z c_i(z,\gamma), \; i=1,2,\dotsc,m \} \cup
\{ \nabla h_j(z), \; j \in \cA(z) \}$ \;\;
is a linearly independent set.}
\end{equation}
The {\em strict complementarity} condition is that
\begin{equation} \label{eq:sc}
\mu_j >0, \;\; \mbox{for all $j \in \cA(z)$}.
\end{equation}
Finally, the {\em second-order sufficient conditions} are
\begin{equation} \label{eq:soc}
d^T W(z,\lambda,\mu,\gamma)d >0 \;\; \mbox{for all $d\in N(z)$ with $d \neq 0$},
\end{equation}
where the subspace $N(z)$ is defined as follows:
\begin{equation} \label{eq:def.N}
N(z) := \left\{ d \, : \parbox{.57\linewidth}{\centering$\nabla_zc_i(z,\gamma)^Td=0, \; i=1,2,\dotsc,m$ \\ and \\
$\nabla h_j(z)^Td=0, \; \forall j \in {\cal A}(z)$} \right\},
\end{equation}
and the matrix $W(z,\lambda,\mu,\gamma)$ is the Hessian of the
Lagrangian function for the problem \eqref{eq:genF}, that is,
\begin{equation} \label{eq:defW}
\begin{aligned}
W(&z,\lambda,\mu,\gamma) \, :=  \\
   & \nabla^2f(z) - \sum_{i=1}^m \lambda_i \nabla^2_{zz} c_i(z,\gamma) - \sum_{j\in \cA(z)} \mu_j \nabla^2 h_j(z).
\end{aligned}\end{equation}

In the neighborhood of a value of $\gamma$ at which all the conditions
above are satisfied, we find that $(z,\lambda,\mu_\cA)$ is an implicit
function of $\gamma$. We find expressions for the derivatives of
$(z,\lambda,\mu_\cA)$ with respect to $\gamma$ by applying the
implicit function theorem (see for example \cite[Theorem~A.2]{NocW06})
to the equality conditions in \eqref{eq:kkt}, which we can formulate
as follows:
\begin{subequations} \label{eq:kktA}
\begin{align}
\nabla f(z) - \sum_{i=1}^m \lambda_i \nabla_z c_i(z,\gamma) & - \sum_{j \in \cA(z)} \mu_j \nabla h_j(z) = 0, \\
c_i(z,\gamma) & =0, \;\; i=1,2,\dotsc,m, \\
h_j(z) & = 0, \;\; j \in \cA(z).
\end{align}
\end{subequations}
Note that this system of equations is square, with $n+m+|\cA|$
equations and unknowns. Moreover, standard analysis of optimality
conditions shows that its square Jacobian matrix is nonsingular, under
the LICQ \eqref{eq:licq} and second-order sufficient \eqref{eq:soc}
conditions. The implicit function theorem now yields the following:
\begin{equation} \label{eq:dz}
\left[ \begin{matrix}
\nabla z(\gamma) \\ \nabla \lambda (\gamma) \\ \nabla \mu_{\cA} (\gamma)
\end{matrix} \right] =
H(z,\lambda,\mu,\gamma)^{-1}
\left[ \begin{matrix}
\sum_{i=1}^m \lambda_i \nabla_{z \gamma} c_i(z,\gamma) \\
-\nabla_{\gamma} c(z,\gamma)^T \\
0
\end{matrix} \right],
\end{equation}
where
\[
  H(z,\lambda,\mu,\gamma) = \left[ \begin{matrix}
W(z,\lambda,\mu,\gamma) & -\nabla_z c(z,\gamma) & -\nabla h_{\cA}(z) \\
\nabla_z c(z,\gamma)^T & 0 & 0 \\
\nabla h_{\cA}(z)^T & 0 & 0
\end{matrix} \right].
\]
We can derive $\nabla \cF(z)$ from $\nabla z (\gamma)$ through the
definition \eqref{eq:genF}, as follows:
\begin{equation} \label{eq:dF}
\nabla \cF(\gamma) =  \nabla z(\gamma)^T \nabla f(z(\gamma)).
\end{equation}

\subsection{Heuristics for Initializing the Power Adjustment Model}
\label{app:safe.vulnerable}

\subsubsection{Determining Safe and Vulnerable Lines}\label{sec:alg.line}

We have noted that a typical grid contains many ``safe'' lines, for
which large changes to the impedance do not affect the ability of the
grid to serve demands. We discuss here a filtering approach to
identify the complementary set of ``vulnerable'' lines, for which
impedance change causes significant disruption of the grid.
(We assume that the number of vulnerable lines is relatively small;
otherwise, the grid has a systemic vulnerability and would be hard to
defend.)

A naive approach for identifying safe and vulnerable lines is
enumeration. For each line $i$, we set $\gamma=\ub{\gamma}e_i$ (where
$e_i$ is a vector of all zeros except for $1$ in the $i$th entry) and
evaluate $\cF_L(\gamma)$ defined by \eqref{eq:feas.rest}. Vulnerable
lines are taken to be those for which $\cF_L(\gamma)>0$. This
enumeration approach is not very effective, in part because of its
cost (it requires solution of $|\cL|$ different power flow problems)
and because it cannot identify {\em combinations} of lines that are
individually ``safe'' but which together create a vulnerability in the
network.  We therefore propose an alternative heuristic called ESL
(for ``eliminating safe lines'').

The motivation of the ESL heuristic is as follow. If a system
operator, instead of an attacker, is asked to increase impedance on
exactly $\kappa$ lines, then he will choose those $\kappa$ lines so as
to {\em minimize} the disruption to the system.  These lines will not
be an attractive choice on the attacker's side, especially when
increase of impedances on these lines does not lead to any load
shedding, so in general we declare these lines to be ``safe.'' On the
other hand, the lines {\em not} selected by the operator correspond to
those that may cause disruption and hence are an attractive target for
attack. We declare such lines to be ``vulnerable.''

In ESL, following the experiment graphed in
Figure~\subref*{fig:30.perturb.safely}, we seek the value of $\kappa$ in
\eqref{eq:model.LS.total} such that a total perturbation of size
$\kappa\ub{\gamma}$ can be distributed to lines with little load shedding.
Lines $i$ for which $\gamma_i\approx\ub{\gamma}$ are declared to be safe. This
process is repeated until relatively few vulnerable lines remain.  The
appropriate value of $\kappa$ can be found by binary search, by solving the
following modification of the problem in \eqref{eq:model.safe}, which
depends on a working set $\cW \subset \cL$ of lines not yet classified as
safe:

\bseq \label{eq:safe.lines}
\begin{alignat}{2}
 \cH_\cW(\kappa,\ub{\gamma}) = \omin{x,y,\gamma} & p^Ty\\
       \ost & F_L(x,y;\gamma) = 0\\
            & e^T \gamma = \kappa\ub{\gamma} \label{eq:safe.lines.perturb}\\
            & \lb{x} \le x \le \ub{x}\\
            & 0 \le y \le \ub{y}\\
            & 0 \le \gamma_i \le \ub{\gamma} &&~~~ i\in\cW \label{eq:safe.line.w}\\
            & \gamma_i = 0                         &&~~~ i\notin\cW.\label{eq:safe.line.not.w}
\end{alignat}\eseq

\begin{algorithm}[t]\footnotesize
\caption{$\mbox{\sc Eliminating Safe Lines (ESL)}$}\label{alg:ESL}
\begin{algorithmic}[1]
\Require
  \Statex $\cL:$ Set of all lines;
  \Statex $\eta \in (0,1)$: Threshold for screening;
\Ensure
  \Statex $\cW$: a set of vulnerable lines;
  \Statex $\cS$: a set of safe lines;
\medskip
\State $\cS\gets \emptyset$;
\Repeat
\State $\cW \gets \cL\backslash\cS$;
\State Define $\kappa^*$ to be the largest value of $\kappa$
for which  $\cH_\cW(\kappa,\ub{\gamma})<\epsilon$;
\State $\cS'\gets\left\{i\, : \, {\gamma_i}/{\ub{\gamma}}\ge \eta,~ i\in\cW\right\}$;
  \Comment \mbox{newly determined ``safe'' lines}
\State $\cS\gets\cS\cup\cS'$;
\Until{$\cS'=\emptyset$}
\end{algorithmic}
\end{algorithm}

The complete procedure is shown in Algorithm~\ref{alg:ESL}.  We start
by putting all lines $\cL$ into the working set $\cW$, then
successively eliminating from $\cW$ those lines $i$ for which the
solution of \eqref{eq:safe.lines} yields $\gamma_i$ within a factor
$\eta$ of the upper bound $\ub{\gamma}$. (We used $\eta=.9$.) The
process is repeated until no new ``safe'' lines are identified. The
lines remaining in $\cW$ are then classified as ``vulnerable.''


\subsubsection{Target Node Selection}\label{sec:alg.node}

\begin{figure}\centering%
  \includepgfplots[width=.85\linewidth]{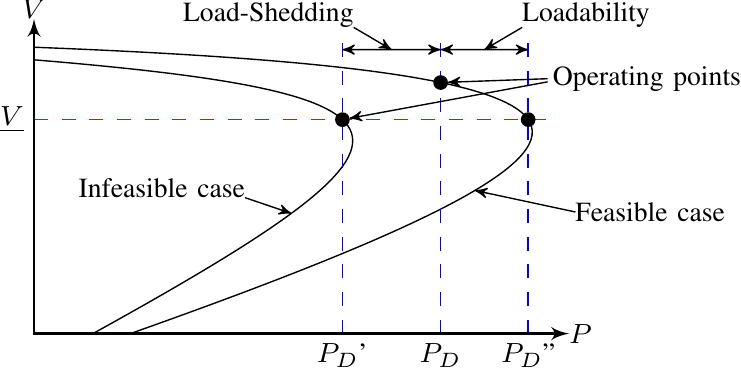}%
  \caption{Loadability and Load Shedding in $PV$-curve}\label{fig:ML.LS}
\end{figure}

Our approach for selecting ``target'' nodes $\cT$ in
Algorithm~\ref{alg:VA.LS} is based on maximum loadability.  The
maximum loadability problem is similar to feasibility restoration in
that it seeks the boundary of the feasible region. However, rather
than starting from an infeasible point (where the nominal loads cannot
be served), it begins from a feasible grid and increases the loads
until demands can no longer be met. The difference is illustrated in
Figure~\ref{fig:ML.LS}, which shows the $PV$-curve for a particular
demand node. When the grid is feasible with demand $P_D$, (right
curve), the demand can be increased to $P_D''$ while retaining
feasibility. The difference $P_D''-P_D$ can be regarded as the maximum
loadability at this node. If the grid is infeasible (left curve) the
demand must be reduced to $P_D'$ before feasibility is recovered.

The formulation for maximum loadability problem can be obtained by
replacing \eqref{eq:feas.rest.power.sigma.p},
\eqref{eq:feas.rest.power.sigma.m}, and \eqref{eq:feas.rest.power.rho}
by the following constraints:
\bseq\label{eq:max.load}\begin{alignat}{2}
  \sigma_i^+ &=~ 0       &&~~~ i\in \cG\label{eq:max.load.sigma.p}\\
  \sigma_i^- &=~ 0       &&~~~ i\in \cG\label{eq:max.load.sigma.m}\\
  \rho_i     &\le~ 0     &&~~~ i\in \cD.\label{eq:max.load.rho}
\end{alignat}\eseq
The first two constraints fix the power generations at their nominal
values, while \eqref{eq:max.load.rho} allows increase (rather than
decrease) of demand at the demand buses.
When the nominal loads and generations are feasible, we
expect the objective to be negative at the solution.

To identify the target nodes, we simply set $\gamma_i=\ub\gamma$ for
all vulnerable lines $i \in \cW$ that are identified by the ESL
procedure, Algorithm~\ref{alg:ESL}, and solve
\eqref{eq:feas.rest.power} for this value of $\gamma$. If a node does
not require any load shedding under this maximal-perturbation setting,
it is unlikely that any attack on the vulnerable lines will lead to
load shedding on this node.  The target nodes are defined to be those
for which load shedding is required, that is, $\rho_i>0$ at the
solution of \eqref{eq:feas.rest.power}.  We denote the set of these
nodes by $\cT$.

\begin{figure}
\centering
 \subfloat[Non-Target Node (No Load Shedding Required)]{%
    \includepgfplots[width=.8\linewidth]{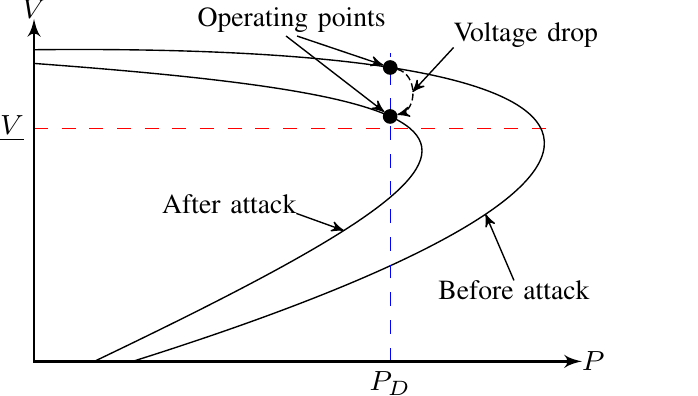}%
    \label{fig:pv.attack.1}}\\
 \subfloat[Target Node (Load Shedding Required)]{%
    \includepgfplots[width=.8\linewidth]{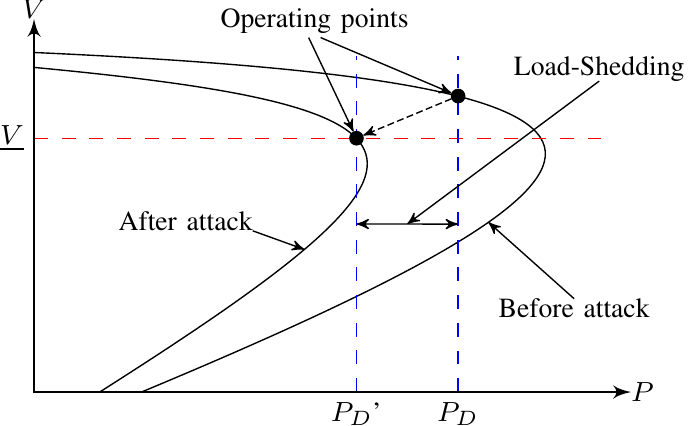}%
    \label{fig:pv.attack.2}}
\caption{Possible Changes of $PV$-curve on a Node After Attack.}\label{fig:pv.attack}
\end{figure}

Figure~\ref{fig:pv.attack} shows target and non-target nodes, and
shows how maximum loadability motivates their classification into
these categories.  The top figure shows a non-target node, for which
is it possible to meet the original demand $P_D$ even after the
maximum-perturbation attack, though the loadability is decreased.
The bottom figure shows that the demand must be reduced to $P_D'$ in
order for the network to remain feasible.  On this node, there is a
chance that an attack on the vulnerable lines will lead to load
shedding.  The target nodes are the nodes affected by the type of
attack we are considering, so that even for a choice of $\gamma$ that
allows the nominal demand on these loads to be served, the change in
maximum loadability may give us some information on the sensitivity of
demand that can be served to the value of $\gamma$.  Since the
objective in \eqref{eq:feas.rest.power} is of weighted $\ell_1$ type,
we expect the number of target nodes to be small.

The target nodes $\cT$ are incorporated into Algorithm~\ref{alg:VA.LS}
by replacing the lower bound \eqref{eq:feas.rest.power.rho} in the
formulation \eqref{eq:feas.rest.power} by a negative quantity for the
nodes in $\cT$, and increasing these bounds toward zero progressively
during the course of Algorithm~\ref{alg:VA.LS}. Additional details are
given in Subsection~\ref{sec:alg.LS.init}.

\IEEEtriggeratref{2}
\IEEEtriggercmd{\enlargethispage{-3.3in}}
\bibliographystyle{IEEEtran}
\bibliography{vulnerability_v2}

\begin{thebibliography}{10}
\providecommand{\url}[1]{#1}
\csname url@samestyle\endcsname
\providecommand{\newblock}{\relax}
\providecommand{\bibinfo}[2]{#2}
\providecommand{\BIBentrySTDinterwordspacing}{\spaceskip=0pt\relax}
\providecommand{\BIBentryALTinterwordstretchfactor}{4}
\providecommand{\BIBentryALTinterwordspacing}{\spaceskip=\fontdimen2\font plus
\BIBentryALTinterwordstretchfactor\fontdimen3\font minus
  \fontdimen4\font\relax}
\providecommand{\BIBforeignlanguage}[2]{{%
\expandafter\ifx\csname l@#1\endcsname\relax
\typeout{** WARNING: IEEEtran.bst: No hyphenation pattern has been}%
\typeout{** loaded for the language `#1'. Using the pattern for}%
\typeout{** the default language instead.}%
\else
\language=\csname l@#1\endcsname
\fi
#2}}
\providecommand{\BIBdecl}{\relax}
\BIBdecl

\bibitem{SalWB04}
J.~Salmeron, K.~Wood, and R.~Baldick, ``{Analysis of electric grid security
  under terrorist threat},'' \emph{IEEE Transactions on Power Systems},
  vol.~19, no.~2, pp. 905--912, May 2004.

\bibitem{MotAG05}
A.~L. Motto, J.~M. Arroyo, and F.~D. Galiana, ``{A Mixed-Integer LP procedure
  for the analysis of electric grid security under disruptive threat},''
  \emph{IEEE Transactions on Power Systems}, vol.~20, no.~3, pp. 1357--1365,
  Aug. 2005.

\bibitem{ArrG05}
J.~M. Arroyo and F.~D. Galiana, ``{On the solution of the bilevel programming
  formulation of the terrorist threat problem},'' \emph{IEEE Transactions on
  Power Systems}, vol.~20, no.~2, pp. 789--797, May 2005.

\bibitem{SalWB09}
J.~Salmeron, K.~Wood, and R.~Baldick, ``{Worst-case interdiction analysis of
  large-scale electric power grids},'' \emph{IEEE Transactions on Power
  Systems}, vol.~24, no.~1, pp. 96--104, Feb. 2009.

\bibitem{DonLL08}
V.~Donde, V.~L\'{o}pez, and B.~C. Lesieutre, ``{Severe multiple contingency
  screening in electric power systems},'' \emph{IEEE Transactions on Power
  Systems}, vol.~23, no.~2, pp. 406--417, May 2008.

\bibitem{PinMD10}
A.~Pinar, J.~Meza, V.~Donde, and B.~C. Lesieutre, ``{Optimization strategies
  for the vulnerability analysis of the electric power grid},'' \emph{SIAM
  Journal on Optimization}, vol.~20, no.~4, pp. 1786--1810, 2010.

\bibitem{Arr10}
J.~Arroyo, ``{Bilevel programming applied to power system vulnerability
  analysis under multiple contingencies},'' \emph{IET Generation, Transmission
  \& Distribution}, vol.~4, no.~2, pp. 178--190, Sep. 2010.

\bibitem{DelAA10}
A.~Delgadillo, J.~M. Arroyo, and N.~Alguacil, ``{Analysis of electric grid
  interdiction with line switching},'' \emph{IEEE Transactions on Power
  Systems}, vol.~25, no.~2, pp. 633--641, May 2010.

\bibitem{ArrF13}
J.~M. Arroyo and F.~J. Fern\'{a}ndez, ``{Application of a genetic algorithm to
  $n-K$ power system security assessment},'' \emph{International Journal of
  Electrical Power \& Energy Systems}, vol.~49, pp. 114--121, Jul. 2013.

\bibitem{BieV10}
D.~Bienstock and A.~Verma, ``{The $N-k$ problem in power grids: new models,
  formulations, and numerical experiments},'' \emph{SIAM Journal on
  Optimization}, vol.~20, no.~5, pp. 2352--2380, 2010.

\bibitem{usc04}
\BIBentryALTinterwordspacing
{U.S.-Canada Power System Outage Task Force}, ``Report on the august 14, 2003
  blackout in the united states and canada: Causes and recommendations,'' 2004.
  [Online]. Available: \url{https://reports.energy.gov}
\BIBentrySTDinterwordspacing

\bibitem{BerV99}
A.~R. Bergen and V.~Vittal, \emph{{Power systems analysis}}, 2nd~ed.\hskip 1em
  plus 0.5em minus 0.4em\relax Prentice Hall, Aug. 1999.

\bibitem{IbaIT80}
K.~Iba, S.~Iwamoto, and Y.~Tamura, ``{A method of finding multiple load-flow
  solutions for actual power systems},'' \emph{Electrical Engineering in
  Japan}, vol. 100, no.~3, pp. 257--264, May 1980.

\bibitem{FraW56}
M.~Frank and P.~Wolfe, ``{An algorithm for quadratic programming},''
  \emph{Naval research logistics quarterly}, vol.~3, no. 1-2, pp. 95--110, Mar.
  1956.

\bibitem{Dunn80a}
J.~C. Dunn, ``Convergence rates for conditional gradient sequences generated by
  implicit step length rules,'' \emph{SIAM Journal on Control and
  Optimization}, vol.~18, no.~5, pp. 473--487, 1980.

\bibitem{AlsS74}
O.~Alsac and B.~Stott, ``{Optimal load flow with steady-state security},''
  \emph{IEEE Transactions on Power Apparatus and Systems}, vol. PAS-93, no.~3,
  pp. 745--751, May 1974.

\bibitem{WacB06}
A.~W\"{a}chter and L.~Biegler, ``{On the implementation of an interior-point
  filter line-search algorithm for large-scale nonlinear programming},''
  \emph{Mathematical Programming}, vol. 106, no.~1, pp. 25--57, Mar. 2006.

\bibitem{ZimMT11}
R.~D. Zimmerman, C.~E. Murillo-S\'{a}nchez, and R.~J. Thomas, ``{MATPOWER:
  Steady-state operations, planning, and analysis tools for power systems
  research and education},'' \emph{IEEE Transactions on Power Systems},
  vol.~26, no.~1, pp. 12--19, Feb. 2011.

\bibitem{NocW06}
J.~Nocedal and S.~J. Wright, \emph{{Numerical Optimization}}, 2nd~ed.\hskip 1em
  plus 0.5em minus 0.4em\relax New York: Springer, 2006.

\end{thebibliography}

\begin{IEEEbiographynophoto}{Taedong Kim}
received the B.S. degree in Computer Science and Engineering from
Seoul National University, Seoul, South Korea in 2007, and M.S. degree
in Computer Sciences from the University of Wisconsin-Madison
(UW-Madison) in 2010. He is currently pursuing the Ph.D degree in
Computer Sciences at UW-Madison.  His research interests lie on
applications of numerical optimization techniques to problems in
sciences and engineering.
\end{IEEEbiographynophoto}
\begin{IEEEbiographynophoto}{Stephen J. Wright}
received the B.Sc. (Hons.) and Ph.D. degrees from the University of
Queensland, Australia, in 1981 and 1984, respectively.

After holding positions at North Carolina State University, Argonne
National Laboratory, and the University of Chicago, he joined the
Computer Sciences Department at the University of Wisconsin-Madison as
a Professor in 2001. His research interests include theory,
algorithms, and applications of computational optimization.

Dr. Wright was Chair of the Mathematical Programming Society from
2007-2010 and served from 2005-2014 on the Board of Trustees of the
Society for Industrial and Applied Mathematics (SIAM).  He has served
on the editorial boards of Mathematical Programming (Series A), SIAM
Review, and the SIAM Journal on Scientific Computing. He has been
editor-in-chief of Mathematical Programming (Series B) and is current
editor-in-chief of the SIAM Journal on Optimization.
\end{IEEEbiographynophoto}
\begin{IEEEbiographynophoto}{Daniel Bienstock}
is a professor at the Departments of Industrial Engineering and
Operations Research and Department of Applied Physics and Applied
Mathematics, Columbia University, where he has been since 1989.  His
research focuses on optimization and computing, with special interest
in power grid modeling and analysis.
\end{IEEEbiographynophoto}
\begin{IEEEbiographynophoto}{Sean Harnett}
is a PhD student at the Department of Applied Physics and Applied
Mathematics, Columbia University.
\end{IEEEbiographynophoto}
\enlargethispage{-3.3in}

\end{document}